\outer\def\theorem#1. {\medbreak\noindent{\bf 
#1.\enspace}\begingroup\sl}
\def\endtheorem{\rm \par
                 \ifdim\lastskip<\medskipamount 
                   \removelastskip \penalty 55 \smallskip
                 \fi
                 \endgroup
                }

\def\pf{\medskip\noindent{\it Proof\/}.\enspace}
\def\qed{\qquad \vrule height6pt width5pt depth0pt}
\def\hat{\widehat{\phantom{m}}\,}
\def\tild{\widetilde{\phantom{m}}\,}
\def\z{{\bf Z}}

\def\n{{\bf N}}

\def\oh{{\cal O}}
\def\a{{\cal A}}
\def\e{{\cal E}}
\def\u{{\cal U}}
\def\g{{\cal G}}
\def\t{{\cal T}}
\def\w{{\cal W}}
\def\k{{\cal K}}
\def\i{{\cal I}}
\def\j{{\cal J}}
\def\tildee{{\widetilde E}}
\def\tildef{\widetilde F}
\def\tildes{\widetilde S}

\def\tilden{\widetilde N}
\def\restrictedto#1{\Big|\lower2pt\hbox{$#1$}}
\def\subsetnoteq{\,\rlap{\raise4true pt\hbox{$\scriptstyle\subset$}}
     \lower3true pt\hbox{$\scriptstyle\not=$}\,}
\def\cstar{$C^*$}
\def\Chi{\raise3true pt\hbox{$\chi$}}    
\def\supp{\hbox{supp}\;}
\def\inv{^{-1}}
\def\ofte{\bigl(t(e)\bigr)}
\def\oftf{\bigl(t(f)\bigr)}
\def\proj{\hbox{Proj}\;}
\def\subdel{_\partial}
\def\implies{\;\Longrightarrow\;}

\def\bib#1{{\bf [#1]}}
\def\anantharaman{1}
\def\anantharamanrenault{2}
\def\archboldspielberg{3}
\def\bprs{4}
\def\brenken{5}
\def\cuntz{6}
\def\cuntzkrieger{7}
\def\drinentomforde{8}
\def\exellaca{9}
\def\flr{10}

\def\higgins{11}
\def\kelley{12}
\def\kprr{13}
\def\kpr{14}
\def\kumjianpask{15}
\def\lacaspielberg{16}
\def\mrw{17}
\def\quigg{18}
\def\quiggsieben{19}
\def\raeburnszymanski{20}
\def\renault{21}
\def\serre{22}

\def\spielberg{23}
\def\szymanskizhang{23}

\centerline{\bf A Functorial Approach to the $C^*$-algebras of a Graph}
\bigskip\footnote{}{1991 {\it Mathematics Subject 
Classification}.  Primary 46L55}
\noindent
\hskip2.5true in Jack Spielberg\par\noindent
\hskip2.5true in Department of Mathematics\par\noindent
\hskip2.5true in Arizona State University\par\noindent
\hskip2.5true in Tempe, AZ  85287-1804\par\noindent
\hskip2.5true in {\it jss@math.la.asu.edu\/}
\bigskip
\centerline{\bf Abstract.}
\smallskip
{\narrower
A functor from the category of directed trees with inclusions to the 
category of commutative \cstar-algebras with injective 
$*$-homomorphisms is constructed.  This is used to define a functor 
from the category of directed graphs with inclusions to the category 
of \cstar-algebras with injective 
$*$-homomorphisms.  The resulting \cstar-algebras are identified as 
Toeplitz graph algebras.  Graph algebras are proved to have 
inductive limit decompositions over any family of subgraphs with 
union equal to the whole graph.  The construction is used to prove 
various structural properties of graph algebras.
\bigskip}
\openup.9\jot
\centerline{\bf Introduction.}
\smallskip
Since the paper of Cuntz and Krieger in 1976, much work has gone into 
elucidating the brief remarks made there regarding the case of 
infinite $0-1$ matrices.  While perhaps the most far-reaching 
solution put forward has been a direct generalization to infinite 
$0-1$ matrices (\bib\exellaca), most of the papers on the subject generalize 
to a class of infinite directed graphs.  (In fact, this is the 
direction indicated in \bib\cuntzkrieger.)  In this paper we propose a new method 
of accomplishing this generalization to infinite graphs, based on a 
construction of a boundary for an arbitrary directed tree.  For a 
locally finite tree, (or more generally, a directed tree in which each 
vertex emits finitely many edges), the obvious space of (directed) 
ends serves as a boundary.  Without the finiteness condition, 
however, the space of ends is not locally compact.  This difficulty 
has been solved in certain cases by ad-hoc methods (e.g., 
\bib{\renault,\kprr,\szymanskizhang}).
More recently these difficulties have been overcome more generally, 
though again by ad-hoc methods 
(\bib{\drinentomforde,\raeburnszymanski}).
Our 
approach is to give a natural solution in the form of a functor from 
the category of directed trees (with inclusion) to the category of 
commutative $C^*$-algebras (with injective $*$-homomorphisms).  
In the case of a locally finite tree,
this functor yields the algebra of continuous 
functions vanishing at infinity on the space of ends.  Given a 
directed graph, the (combinatorial) universal cover is a bundle of 
directed trees, on which acts the (combinatorial) fundamental 
groupoid.  Applying our boundary functor, we may let the groupoid act 
on the resulting $C^*$-algebra.  We prove that this yields the same 
graph algebra (up to strong Morita equivalence) as has been defined 
by others. 
 The chief virtue of our approach seems to us to be its 
naturality.  We obtain a morphism of graph algebras from an inclusion 
of directed graphs.  This leads us to introduce a family of Toeplitz 
graph algebras associated to a directed graph. We immediately obtain 
an inductive limit decomposition of a graph algebra over any directed 
family of subgraphs.  We use this theorem in \bib\spielberg\ to establish 
semiprojectivity for classifiable simple purely infinite 
$C^*$-algebras having finitely generated $K$-theory and torsion-free 
$K_1$.  Moreover, we mention that all of the structure of the graph 
algebras is obtained without any restrictions on the size or structure 
of the graph.  (A different method for obtaining the direct limit 
decomposition over the finite subgraphs is given in 
\bib{\raeburnszymanski}.)

The outline of the paper is as follows.  We finish the introduction 
with the basic notation we will use.  We find it convenient to 
separate the {\it graphical} properties of an edge from the {\it 
groupoid} properties.  Thus we adopt the terminology of {\it origin} 
and {\it terminus}, as in \bib\serre, leaving {\it source} and {\it 
range} exclusively as groupoid terms.  In section 1 we introduce our 
boundary functor for directed trees.  We include a classification of 
the open subsets of the boundary, as this is needed to classify the 
ideals in the graph algebras.  In section 2 we construct the groupoids 
and $C^*$-algebras obtained from a directed graph by following the 
above approach.  We prove the standard presentation by generators and 
relations, and define the fundamental cocycle.  This allows us to 
exhibit the AF core as the $C^*$-algebra of an equivalence relation, 
and prove nuclearity for all graph algebras.  We also prove the 
inductive limit decomposition described above.  In section 3 we prove 
the basic structure theorems regarding simplicity and pure 
infiniteness.  These properties follow quite easily from the theorems 
already established.

The idea of using a Boolean ring of subsets of the vertices of a(n 
undirected) tree 
in order to define its boundary was discovered jointly with Marcelo 
Laca, when the author was visiting the Mathematics Department at the 
University of Newcastle.  He wishes to thank the members of the 
Department, and Marcelo Laca and Iain Raeburn in particular, for their 
wonderful generosity and hospitality during that stay.
\medskip

Following \bib\serre\ we let a {\it graph}, $E$, consist of a set $E^0$ of {\it 
vertices}, a set $E^1$ of {\it edges}, maps $o$, $t:E^1\to E^0$
({\it origin\/} and {\it terminus\/}), and 
a map $e\mapsto \overline e$ on $E^1$ ({\it reversal\/}), satisfying:
$$\eqalign{\overline{\overline e}&=e\cr\overline e&\not=e\cr 
o(\overline e)&=t(e).\cr}$$
A {\it directed graph} consists of a graph $E$ 
together with a subset $E^1_+\subseteq E^1$ containing exactly one 
edge from each pair $\{e,\overline e\}$.
(When drawing $E$ it is usual to give only one segment to represent 
$\{e,\overline e\}$.
If $E$ is directed, the segment representing $\{e,\overline e\}$ is 
given an arrowhead pointing toward $t(|e|)$, where $|e|$ is the 
element of $\{e,\overline e\}\cap E^1_+$.)  
A {\it path} in $E$ is a 
(finite or infinite) sequence of edges $e_1e_2\cdots$ such that 
$t(e_{i-1})=o(e_i)$ and $e_i\not=\overline{e_{i-1}}$
 for all $i>1$.  We write $o(e_1e_2\cdots) 
=o(e_1)$, $t(e_1\cdots e_n)=t(e_n)$, and $\overline{e_1 
\cdots e_n} =\overline{e_n}\cdots\overline{e_1}$.
If $p$ is a path in $E$ we 
let $\ell(p)$ denote the {\it length} of $p$ ($0\le\ell(p)\le\infty$).
  If $E$ is a directed graph 
we call a path $e_1e_2\cdots$ {\it directed} if $e_i\in E^1_+$ 
for all $i$.  (A path of length zero consists of a single vertex, and 
is considered to be directed.)
We let $E^n$ (respectively $E^n_+$) denote the paths 
(respectively directed paths) of length $n$, for $0\le 
n\le\infty$, $E^*=\bigcup_{0\le n<\infty}E^n$, $E_+^*=\bigcup_{0\le n
<\infty}E_+^n$, and $E^{**}=E^*\cup E^\infty$,
$E^{**}_+=E^*_+\cup E^\infty_+$.
The graph $E$ is a {\it tree} if for all $u,v\in 
E^0$ there exists a unique path $p\in E^*$ with $o(p)=u$ and $t(p)=v$.
A graph is a {\it forest} if it is a disjoint union of trees.  More 
precisely, the graph $E$ is a forest if for all $u,v\in 
E^0$ there exists at most one path $p\in E^*$ with $o(p)=u$ and $t(p)=v$.

\smallskip
Let $E$ be a directed graph.  For $v\in E^0$ let
$$\eqalignno{
V(v)\equiv V_E(v)&=\bigl\{t(p)\bigm|p\in E^*_+,\;o(p)=v\bigr\},\cr
\Delta_1(v)\equiv\Delta_{1,E}(v)
         &=\bigl\{e\in E^1_+\bigm|o(e)=v\bigr\}\cr
\Omega\equiv\Omega(E)&=\bigl\{u\in 
E^0\bigm|\Delta_1(u)=\emptyset\bigr\}\cr
\Sigma\equiv\Sigma(E)&=\bigl\{u\in E^0\bigm| 0<\#\Delta_1(u) 
<\infty\bigr\},\cr
\noalign{\hbox{and\ for\ $F\subseteq\Delta_1(v)$\ a\ finite\ subset,}}
V(v;F)\equiv V_E(v;F)&=\{v\}\cup\bigcup_{e\in\Delta_1(v)\setminus
 F}V\ofte.\cr
		 }$$
Thus $V(v)$ is the set of vertices that can be found by following a 
finite directed path from $v$, $V(v;F)$ is the set of vertices 
that can be found by following a directed path from $v$ that does not 
begin with an edge from the finite set $F$, and $\Omega$ is the set of 
vertices often termed {\it sinks}.  If  $e\in 
\Delta_1(v)$, we will also write $V(v;e)$ for $V(v;\{e\})$.
We will also use the following
notations for paths, respectively directed paths, ending at a sink:
$$\eqalign{\Omega^*&=t\inv(\Omega)\cap E^*\cr
\Omega^*_+&=\Omega^*\cap E^*_+.\cr}$$

\bigskip
\centerline{\bf 1.  The Boundary of a Directed Tree.}
\smallskip\noindent
{\bf Remark  1.1.}\enspace Suppose that $T$ is a directed tree, and $v\in T^0$.
  If $e$, $f\in\Delta_1(v)$ are distinct,
then $V\ofte\cap 
V\oftf=\emptyset$.
\smallskip
It follows that for 
 a directed tree $T$, we may also write
$$V_T(v;F) = V_T(v)\setminus\bigcup_{e\in F}V_T\ofte.$$

\smallskip
For $T$  a directed tree, we define
$$\eqalign{
\e\equiv\e(T)&=\bigl\{V(v;F)\bigm|v\in T^0,\;
F\subseteq\Delta_1(v) \hbox{\ 
finite}\bigr\}\cup\bigl\{\emptyset\bigr\},\cr
\a\equiv\a(T)&=\bigr\{\cup_{i=1}^kB_i\bigm|B_1,\ldots,B_k\in\e \hbox{\ 
are\ disjoint},\;k\ge1\bigr\}.\cr}$$

\theorem Lemma  1.2. If $B$, $C\in\e$, then $B\cap C$, $B\setminus C\in\a$.
\endtheorem

\pf Let $B=V(u;F)$, $C=V(v;G)$, where $F=\{f_1,\ldots,f_k\}$, 
$G=\{g_1,\ldots,g_\ell\}$.  Suppose first that $u=v$.  Then $B\cap 
C=V(u;F\cup G)$ and $B\setminus C=\bigcup_{e\in G\setminus F}V\ofte$; both 
are in $\a$ (the second by  Remark 1.1).

\smallskip
Now suppose that $u\not= v$.  We distinguish three cases.
\smallskip\noindent{\it case (i):\/} $v\in V(u)$.  Let 
$e_1\cdots e_n$ be the (unique directed) path from $u$ to $v$.  If 
$e_1\in F$, then $V\bigl(t(e_1)\bigr)\cap B=\emptyset$.  Since 
$C\subseteq V\bigl(t(e_1)\bigr)$, it follows that $B\cap C=\emptyset$, 
$B\setminus C=B$, and $C\setminus B=C$.
If $e_1\not\in F$, then $C\subseteq 
V(v)\subseteq B$, so $B\cap C=C$.  Thus $C\setminus B=\emptyset$, and
$$B\setminus C = V\bigl(u;F\cup\{e_1\}\bigr) \cup 
V\bigl(t(e_1);e_2\bigr)\cup \cdots \cup 
V\bigl(t(e_{n-1});e_n\bigr) \cup \bigcup_{e\in G}V\ofte,$$
a disjoint union of sets from $\e$.

\smallskip\noindent{\it case (ii):\/} $u\in V(v)$.  This case is 
analogous to case (i).

\smallskip\noindent{\it case (iii):\/} $u\not\in V(v)$ and $v\not\in 
V(u)$.  If $V(u)\cap V(v)=\emptyset$, then $B\cap C=\emptyset$ and 
$B\setminus C=B$.  Otherwise, the path from $u$ to $v$ has the form 
$p\overline q$, where $p$, $q\in E^*_+$.  Then $B\cap 
C=V\bigl(t(p)\bigr)$, and $B\setminus C=B\setminus 
V\bigl(t(p)\bigr)$, which falls under case (i).\qed

\smallskip\noindent
\theorem Lemma 1.3.  $\a$ is a (Boolean) ring of sets.\endtheorem
\pf It is a standard fact
that if $\e$ is a collection of sets containing 
$\emptyset$, $\a$ is the collection of finite disjoint unions of sets 
from $\e$, and if for any $B$, $C\in\e$ we have $B\cap C$, $B\setminus 
C\in\a$, then $\a$ is a ring of sets.  \qed

\smallskip
Before continuing, we need some more properties of sets in the 
ring $\a(T)$.
In the following lemma, the notation $V(x;y)$ will 
mean that $x\in T^0$ and that $y$ is a finite subset of $\Delta_1(x)$.

\theorem Lemma 1.4. Let $T$ be a directed tree.
\item{(i)} $V(u;F)\subseteq V(v;G)$ if and only if one of the 
following two conditions holds:
\itemitem{(a)} $u=v$ and $F\supseteq G$, or
\itemitem{(b)} $u\not= v$, and if $p=e_1e_2\cdots\in T^*$ is the 
unique path with $o(p)=v$ and $t(p)=u$, then $p\in T^*_+$ and 
$e_1\not\in G$.
\item{(ii)} If $V(u;F)\subseteq V(v;G)$, then $v\in V(u;F)\iff u=v$.
\item{(iii)} If $V(v;G)=\bigcup_{i=1}^kV(u_i;F_i)$ is a disjoint 
union, then
$$\displaylines{
\ \hbox{(a)}\; u_i=v \hbox{\ for\ some\ $i$,\
say\ }i=1, \hfill\cr
\  \hbox{(b)}\; F_1\supseteq G,\hfill\cr
\ \hbox{(c)}\; \bigcup_{i=2}^k V(u_i;F_i)= \bigcup_{e\in 
F_1\setminus G} V\ofte.\hfill\cr}$$
\endtheorem

\pf (i) The {\it if} is clear.  For the {\it only if},
part (a) is clear.  Suppose that $u\not=v$.
Then part (b) follows from the fact that $u\in V(v;G)$.
\smallskip\noindent
(ii) Let $p$ be as in (i).  Then $p\in T^*_+$ by (i).  Now 
$v\in V(u;F)\iff \overline p\in T^*_+ \iff u=v$.
\smallskip\noindent
(iii) Part (a) follows from (ii).  Then part (b) follows from (ia). 
Finally, for part (c) note that
$$\eqalign{
\bigcup_2^k V(u_i;F_i)&= V(v;G)\setminus V(v;F_1)\cr
&= \bigcup_{e\in F_1\setminus G}V\ofte.\qed\cr}$$

\smallskip
In the next result we charactize Boolean ring homomorphisms from 
$\a(T)$.

\theorem Lemma 1.5.  Let $\Theta$ be a Boolean ring, and for $u\in 
T^0$ let $\theta_u\in\Theta$ satisfy
\smallskip
\item{(i)} $e\in\Delta_1(v)\implies\theta_{t(e)}\subseteq\theta_v$
\item{(ii)} $e_1,e_2\in\Delta_1(v)$ and $e_1\not=e_2\implies
\theta_{t(e_1)}\cap\theta_{t(e_2)}=\emptyset$.
\smallskip\noindent
Then there exists a unique Boolean ring homomorphism $\mu:\a(T) \to 
\Theta$ such that 
$\mu\bigl(V(u)\bigr)$ $=\theta_u$, for $u\in T^0$.
\endtheorem

\pf We claim that the assignment $V(u;F)\mapsto \theta_u \setminus 
\bigcup\{\theta_{t(e)}\mid e\in F\}$ is well-defined.  First, if 
$V(u;F)=V(v;G)$, then by Lemma 1.4 (i) and (ii) we have $u=v$ and 
$F=G$.  Suppose inductively that $\mu(B)=\cup_i\mu(B_i)$ whenever 
$B,B_i\in\e$ and $B=\cup_iB_i$ is a disjoint union of fewer than $n$ 
sets, and 
let $V(v;G)=\cup_{i=1}^nV(u_i;F_i)$ be a disjoint union.  
By Lemma 1.4, possibly after reordering the $\{u_i\}$, we have
$$\displaylines{
\quad u_1=v\hfill\cr
\quad V(v;G) = V(v;F_1)\cup\bigcup_{e\in F_1
\setminus G}V\ofte,\hbox{\ a\ disjoint\ union,}\hfill\cr
\quad V\ofte=\bigcup\bigl\{V(u_i;F_i)\bigm|u_i\in 
V\ofte\bigr\},\hbox{\ for\ }e\in F_1\setminus G.\hfill\cr}$$
Since $G\subseteq F_1$, we have
$$\theta_v\setminus\bigcup_{e\in G}\theta_{t(e)}=\bigl(
\theta_v\setminus\bigcup_{e\in F_1}\theta_{t(e)}\bigr) \cup \bigcup_{f\in 
F_1\setminus G}\theta_{t(f)}.$$
For $e\in F_1\setminus G$, $\#\{i\mid u_i\in V_{T_1}\ofte\}<n$.  By the 
inductive hypothesis, we have for each $e\in F_1\setminus G$,
$$\theta_{t(e)} = \bigcup_{u_i\in V\ofte}\mu\bigl(V(u_i;F_i)\bigr).$$
Thus
$$\eqalign{
\mu\bigl(V(v;G)\bigr)
&=\mu\bigl(V(v;F_1)\bigr)\cup 
\bigcup_{e\in F_1\setminus G}\theta_{t(e)}\cr
&=\mu\bigl(V(v;F_1)\bigr)\cup 
\bigcup_{e\in F_1\setminus G}\bigcup_{u_i\in V\ofte}
\mu\bigl(V(u_i;F_i)\bigr)\cr
&= \bigcup_{i=1}^n \mu\bigl(V(u_i;F_i)\bigr).\cr}$$
This shows that $\mu$ is well-defined on elements of $\e$.
It follows that $\mu$ is well-defined on $\a$ by $\mu(\cup_iB_i) = 
\cup_i\mu(B_i)$ for finite disjoint collections 
$\{B_i\}\subseteq\e$.  By following the proof of Lemma 1.2, it is now 
easy to see that $\mu(B\cap C)=\mu(B)\cap\mu(C)$ and $\mu(B\setminus 
C)=\mu(B)\setminus\mu(C)$ for any $B,C\in\e$, and hence for 
$B,C\in\a$.  Thus $\mu$ is a homomorphism of Boolean rings.  
Uniqueness follows from the fact that $\a$ is generated by 
$\{V(u)\mid u\in T^0\}$. \qed

\smallskip\noindent
{\bf Remark 1.6.}\enspace  With $\mu$ as in Lemma 1.5, $\ker(\mu)$ is the set of 
finite unions of sets of the form $V(u;F)$ for which 
$\theta_u=\bigcup_{e\in F}\theta_{t(e)}$. 

\smallskip
We now let $A\equiv A(T)\subseteq\ell^\infty(T^0)$ be the closed linear 
span of $\bigl\{\Chi_B\bigm|B\in\a\bigr\}$.  Then $A$ is a 
commutative \cstar-algebra, and $\a$ can be identified with
the Boolean ring of projections in $A$.  We note that there is a 
bijective correspondence between the $*$-representations of $A$ into 
a \cstar-algebra $L$, and the ring homomorphisms of $\a$ into the ring 
of projections in a commutative \cstar-subalgebra of $L$.  (The 
correspondence is obtained as follows:  if $\phi:A\to L$ is a 
$*$-homomorphism, then define $\mu:\a\to Proj\bigl(\phi(A)\bigr)$ by 
$\mu(B)=\phi(\Chi_B)$.  Moreover, $\ker \phi=\overline{\hbox{span}}
\{\Chi_B\mid \mu(B)=0\}$.)
We will identify $\widehat A$ with the 
set of $\a$-ultrafilters:  for $\omega\in\widehat A$ we have the 
$\a$-ultrafilter $\{B\in\a\mid\langle\omega,\Chi_B\rangle=1\}$.
(It is elementary to show that this coincides with the Stone space 
of $\a$ (\bib\kelley, 5.S).)
Of course the fixed $\a$-ultrafilters correspond to the elements of 
$T^0$:  for $v\in T^0$, $\u_v=\{B\in\a\mid v\in B\}$.  We will 
characterize the free $\a$-ultrafilters.

\smallskip
For $p=e_1e_2\cdots\in T^\infty_+$ let $\u_p=\bigl\{B\in\a\bigm| 
V\bigl(o(e_j)\bigr)\subseteq B \hbox{\ for\ some\ }j\bigr\}$.  It is 
clear that $\u_p$ is an $\a$-filter.

\theorem Lemma 1.7. $\u_p$ is an $\a$-ultrafilter.\endtheorem

\pf Let $B\in\a$, $B\not=\emptyset$, $B\not\in\u_p$.  
 Write $B=\cup_1^k B_i$ with 
$\{B_i\}\subseteq\e$ disjoint.  Then $B_i\not\in\u_p$ for each $i$.  
Let $B_i=V(v_i;F_i)$, and let $q_i\in T^*$ be the geodesic from $v_i$ 
to $p$ (cf. \bib\serre, 6.3 Lemma 9).  Since $B_i\not\in\u_p$,
then either $q_i\not\in T^*_+$ or the first edge of $q_i$ is in $F_i$.  
In either case, it follows that if $t(q_i)=o(e_{n_i})$, then $B_i\cap 
V\bigl(t(e_{n_i})\bigr)=\emptyset$.  Let $n=\max(n_1,\ldots,n_k)$.  
Then $V\bigl(o(e_n)\bigr)\in\u_p$ and $B\cap V\bigl(o(e_n)\bigr) 
=\emptyset$.  Hence $\u_p$ is a maximal $\a$-filter.\qed

\smallskip\noindent
{\bf Definition 1.8.}\enspace  (cf. \bib\kumjianpask, section 4)
  If $p=e_1e_2\cdots$ and $q=f_1f_2\cdots\in 
T^\infty_+$ we say that $p\sim q$ if there exist $j$ and $k$ so that 
$e_{j+r}=f_{k+r}$ for $r\ge0$.

\smallskip
It is clear that $p\sim q$ if and only if $\u_p=\u_q$.  We let 
$\u_{[p]}$ denote the common $\a$-ultrafilter corresponding to the 
infinite paths in $[p]\in T^\infty_+/\sim$.

\theorem Lemma 1.9.  $\widehat A=T^0\cup(T^\infty_+/\sim)$.\endtheorem

\pf Let $\u$ be a free $\a$-ultrafilter.  Let $V(v)\in\u$, for some 
$v\in T^0$.  We claim that there is $e\in E^1_+$ with $v=o(e)$ and 
$V\bigl(t(e)\bigr)\in\u$.  For suppose not.  Then for any finite set 
$F\subseteq\Delta_1(v)$, $V(v;F)\in\u$ (we have used the fact that if 
a finite union of disjoint sets is contained in an $\a$-ultrafilter, 
then one of the sets must be in the $\a$-ultrafilter).  But then 
$B\in\u$ whenever $B\in\a$ and $v\in B$.  Thus $\u\supseteq\u_v$, and 
hence $\u=\u_v$ is a not a free $\a$-ultrafilter.  This establishes 
the claim.  Since $\u\not=\emptyset$ there exists $v\in T_0$ with 
$V(v)\in\u$.  Then we may use the claim to find $p=e_1e_2\cdots
\in T^\infty_+$ such that $V\bigl(t(e_j)\bigr)\in\u$ for all $j$.  
Then $\u\supseteq\u_{[p]}$.  Since $\u_{[p]}$ is an 
$\a$-ultrafilter, $\u=\u_{[p]}$.\qed

\smallskip\noindent
{\bf Remark 1.10.}\enspace   We will extend the definition of the equivalence 
$\sim$ from $T^\infty_+$ to $T^\infty_+\cup T^*_+$:  for $p$, 
$q\in T^*_+$ we say that $p\sim q$ if $t(p)=t(q)$.  
For example, while it is 
true that $\Omega^*_+/\sim\equiv\Omega$, it will prove useful to have 
the classes represented by paths.

\smallskip
We next define the boundary of a directed tree.  Our definition 
coincides with the usual space of ends in the case where $\Delta_1(u)$ 
is finite for all vertices $u$ (that is, in the case of a {\it 
row-finite\/} tree).
\smallskip\noindent
{\bf Definition 1.11.}\enspace   The {\it boundary} of the directed tree $T$ is the 
closure in $\widehat A$ of the  set of equivalence classes 
of infinite 
directed paths and finite directed paths ending in $\Omega$:
$$\partial T=\overline{(T^\infty_+\cup\Omega^*_+)/\sim}\;.$$

\noindent
{\bf Remark 1.12.}\enspace  For $B\in\a$, we will let $[B]$ denote the support of 
$\Chi_B$ in $\widehat A$ (that is, the set of $\a$-ultrafilters 
containing $B$).  We will let $[B]\subdel=[B]\cap\partial T$.

\smallskip\noindent
{\bf Remark 1.13.}\enspace  For $v\in T^0$, a neighborhood base at $v$ in $\widehat A$ 
is given by $\bigl\{[V(v;F)]\bigm|F\subseteq\Delta_1(v)\hbox{\ 
finite}\bigr\}$.  For $p=e_1e_2\cdots\in T^\infty_+$, a 
neighborhood base at $[p]$ in $\widehat A$ is given by 
$\bigl\{\bigl[V\bigl(o(e_j)\bigr)\bigr]\bigm|j=1,2,\ldots\bigr\}$.

\theorem Lemma 1.14.  Let $v\in T^0$.  Then $v\in\partial T$ iff 
$\Delta_1(v)$ is infinite or empty.\endtheorem

\pf 
If $\Delta_1(v)$ is empty, then 
$v\in \Omega\subseteq\partial T$.  If $\Delta_1(v)$ is infinite, then 
for every finite $F\subseteq\Delta_1(v)$, there is $p=e_1e_2\cdots
\in T^\infty_+ 
\cup\Omega^*_+$ with $e_1\notin F$.  Then 
$[p]\in\bigl[V(v;F)\bigr]$.  Hence 
$v\in\overline{(T^\infty_+\cup\Omega^*_+)/\sim}$.  If $\Delta_1(v)$ is 
nonempty and finite, then $\{v\}=V\bigl(v;\Delta_1(v)\bigr)$.  Hence 
$v\notin\partial T$.  \qed

\smallskip\noindent
{\bf Remark 1.15.}\enspace   $\Sigma(T) = A(T)\hat\setminus\partial T$,
and is a discrete clopen subset of $A(T)\hat$.

\smallskip\noindent
{\bf Definition 1.16.}\enspace  For $S\subseteq\Sigma$ we let 
$$\eqalign{
A(T,S) &= A(T)/C_0(S)\cong C_0\bigl(\partial T\cup(\Sigma\setminus 
S)\bigr)\cr
\partial(T,S)&=A(T,S)\hat\cong S^c.\cr
}$$
Thus $A(T,\emptyset)=A$ and $A(T,\Sigma)=C_0(\partial T)$; 
$\partial(T,\emptyset)=\widehat A$, and $\partial(T,\Sigma)=\partial 
T$.

\theorem
Lemma 1.17.  Let $S\subseteq\Sigma(T)$, let $\pi:A(T)\to A(T,S)$ be 
the quotient map, and let $\theta:\a(T)\to\proj\bigl(A(T,S)\bigr)$ be 
the corresponding Boolean ring homomorphism.  Then $\ker\theta$ equals 
the collection of finite subsets of $S$.\endtheorem

\pf We have $\Sigma(T)\subseteq\a(T)$ via the singleton sets 
$\bigl\{V\bigl(u;\Delta_1(u)\bigr)\bigm|u\in\Sigma(T)\bigr\}$.  The 
finite unions of these sets are precisely the elements of $\a(T)$ 
sent to $\emptyset$ by $\theta$. \qed

\smallskip
Next we discuss the behavior of the algebras $A(T,S)$ under 
mappings of trees.
By a {\it morphism of directed graphs} we mean a pair of maps from 
the set of vertices (edges) of the first
graph to the set of 
vertices (edges) of the second graph
that intertwine $o$, $t$, and 
$\overline{(\cdot)}$, and map directed edges to directed edges.  It 
is a {\it monomorphism} if it is injective on the sets of vertices and 
edges.

\theorem Lemma 1.18.  Let $\alpha:T_1\to T_2$ be a monomorphism of directed 
trees.  Then
$$\eqalign{
V_{T_1}(v;F)&\mapsto \alpha_*\bigl(V_{T_1}(v;F)\bigr)\equiv
 V_{T_2}\bigl(\alpha(v);\alpha(F)\bigr)\cr
\cup_i B_i&\mapsto\cup_i\alpha_*(B_i)\hbox{\ 
($\{B_i\}\subseteq\e(T_1)$\ finite\ 
disjoint)}\cr}$$
defines a Boolean ring monomorphism $\alpha_*:\a(T_1)\to\a(T_2)$.
\endtheorem

\pf We will check the conditions of Lemma 1.5.
Let $e\in\Delta_{1,T_1}(v)$.  Then 
$\alpha(e)\in\Delta_{1,T_2}\bigl(\alpha(v)\bigr)$, so 
$V_{T_2}\bigl(\alpha\ofte\bigr)\subseteq V_{T_2}\bigl(\alpha(v)\bigr)$.  
Now let $e_1\not=e_2\in\Delta_{1,T_1}(v)$.  Since 
$\alpha$ is 
one-to-one on edges, $\alpha(e_1)\ne
\alpha(e_2)\in 
\Delta_{1,T_2}\bigl(\alpha(v)\bigr)$. Thus, since $T_2$ is a tree,
$V_{T_2}\Bigl(t\bigl(\alpha(e_1)\bigr)\Bigr) \cap 
V_{T_2}\Bigl(t\bigl(\alpha(e_2)\bigr)\Bigr) 
=\emptyset$.  By Lemma 1.5 there is a Boolean ring homomorphism 
$\alpha_*:\a(T_1)\to\a(T_2)$ such that 
$\alpha_*\bigl(V_{T_1}(v)\bigr) = V_{T_2}\bigl(\alpha(v)\bigr)$, for 
$v\in T_1^0$.  Then
$$\eqalign{\alpha_*\bigl(V_{T_1}(v;F)\bigr)&=
\alpha_*\bigl(V_{T_1}(v)\setminus\bigcup_{e\in S}V_{T_1}\ofte\bigr)\cr
&=V_{T_2}\bigl(\alpha(v)\bigr)\setminus\bigcup_{e\in S}V_{T_2}\bigl( 
\alpha\ofte\bigr)\cr
&=V_{T_2}\bigl(\alpha(v)\bigr)\setminus\bigcup_{f\in 
\alpha(S)}V_{T_2}\oftf\cr
&=V_{T_2}\bigl(\alpha(v);\alpha(F)\bigr).\cr}$$
 Finally
$\alpha_*$ is a monomorphism since $B$ nonempty implies $\alpha_*(B)$ 
nonempty.\qed

\smallskip\noindent
{\bf Remark 1.19.}\enspace 
If $\alpha:T_1\to T_2$ is a monomorphism of directed trees, then the 
ring monomorphism $\alpha_*:\a(T_1)\to\a(T_2)$ induces an 
injective $*$-homomorphism $\alpha_*:A(T_1)\to A(T_2)$.  

\theorem Lemma 1.20. Let $\alpha:T_1\to T_2$ be a monomorphism of directed 
trees, and let $S_2\subseteq\Sigma(T_2)$.  Let 
$$S_1=\bigl\{v\in \Sigma(T_1)\bigm| \alpha(v)\in S_2\hbox{\ and\ } 
\alpha\bigl(\Delta_{1,T_1}(v)\bigr)
=\Delta_{1,T_2}\bigl(\alpha(v)\bigr)
\bigr\}.$$
Then $\ker(\pi_2 \circ \alpha_*) = 
C_0(S_1)$, so that the following diagram commutes, has surjective columns 
(given by restrictions), and has injective rows:
$$\matrix{A(T_1)&\buildrel\alpha_*\over\longrightarrow&A(T_2)\cr\cr
\Bigm\downarrow&&\Bigm\downarrow \pi_2\cr\cr
A(T_1,S_1)&\longrightarrow&A(T_2,S_2).\cr}$$
\endtheorem

\pf Let $\theta:\a(T_1)\to\proj\bigl(A(T_2,S_2)\bigr)$ be the Boolean 
ring homomorphism associated to $\pi_2 \circ \alpha_*$.  Then 
$V_{T_1}(v;F)\in\ker\theta$ if and only if 
$$\theta\bigl(V_{T_1}(v)\bigr)=\bigcup_{e\in F}\theta\bigl( 
V_{T_1}\ofte\bigr).$$
Since $\alpha_*$ is a monomorphism, this happens if and only if the 
characteristic function of 
$\bigl[V\bigl(\alpha(v);\alpha(F)\bigr)\bigr]$ is in $\ker\pi_2$.  By 
Lemma 1.17, this happens if and only if $\alpha(v)\in S_2$, 
$F=\Delta_{1,T_1}(v)$, and $\alpha\bigl(\Delta_{1,T_1}(v)\bigr) = 
\Delta_{1,T_2}\bigl(\alpha(v)\bigr)$.  \qed

\smallskip\noindent
{\bf Remark 1.21.}\enspace  If $\alpha:T_1\to T_2$ and $\beta:T_2\to T_3$ are 
monomorphisms of directed trees, then it is clear that 
$\beta_*\circ\alpha_*=(\beta\circ\alpha)_*$ at the level of Boolean 
rings, and hence also as $*$-homomorphisms.  Let 
$S_3\subseteq\Sigma(T_3)$.  
Applying Lemma 1.20 to $\beta$ and $S_3$ yields $S_2\subseteq\Sigma(T_2)$.  
Applying Lemma 1.20 to $\alpha$ and $S_2$ yields 
$S_1\subseteq\Sigma(T_1)$.  It is easily verified that 
an application of Lemma 1.20 to $\beta\circ\alpha$ 
and $S_3$ would give the same subset $S_1$ of $\Sigma(T_1)$.  Hence 
the concatenation of commuting squares obtained from Lemma 1.20 gives 
the commuting square for a composition.

\theorem Corollary 1.22. Let $\pi:A(T_2)\to C_0(\partial T_2)$ be the 
restriction mapping. Then $\ker(\pi\circ\alpha_*)=C_0(S)$, where
$$S=\bigl\{v\in \Sigma(T_1)\bigm| \alpha\bigl(\Delta_{1,T_1}(v)\bigr)
=\Delta_{1,T_2}\bigl(\alpha(v)\bigr)\bigr\}.\qed$$\endtheorem

We next characterize the open subsets of the boundary of a directed 
tree.  This will be crucial for determining the ideals in graph 
$C^*$-algebras later on.  We will eventually need the 
generalizations of these 
notions to directed graphs.  Therefore, we state the definitions now 
in that generality.

\smallskip\noindent
{\bf Definition 1.23.}\enspace Let $E$ be a directed graph. An {\it 
invariant} of $E$ is a pair $(N,F)$ consisting of a subset 
$N\subseteq E^0$ and a family $\{F_u\mid u\in N\}$, where 
$F_u\subseteq \Delta_1(u)$ is a finite set, satisfying the following 
conditions.
\smallskip
\item{(i)} If $u\in N$ and $\Delta_1(u)$ is finite, then 
$F_u=\emptyset$.
\smallskip
\item{(ii)} Let $u\in N$ and $e\in\Delta_1(u)$.  Then
\smallskip
\itemitem{(a)} If $e\notin F_u$, then $t(e)\in N$ and 
$F_{t(e)}=\emptyset$.
\smallskip
\itemitem{(b)} If $e\in F_u$ and $t(e)\in N$, then 
$F_{t(e)}\ne\emptyset$.
\smallskip
\item{(iii)} Let $u\in E^0$ with $\Delta_1(u)$ finite and nonempty.  If 
$t\bigl(\Delta_1(u)\bigr)\subseteq N$, and $F_{t(e)}=\emptyset$ 
for all $e\in \Delta_1(u)$, then $u\in N$.
\medskip
\noindent
We let $\i\equiv\i(E)$ denote the set of invariants of $E$.
We remark that the notions of {\it hereditary} and {\it saturated} 
subsets of $E^0$ appearing in other work are being replaced by 
conditions (i) and (ii)(a), respectively (ii)(b) and (iii), in the 
context of general directed graphs.
We
define a partial ordering in $\i$ by
$(N_1,F_1)\le(N_2,F_2)$ if $N_1\subseteq N_2$ and $F_{1,u}\supseteq 
F_{2,u}$ for each $u\in N_1$.
\smallskip\noindent
{\bf Example 1.24.}\enspace Let $E^0 = \z$ and let $E$ have infinitely 
many directed edges from $n$ to $n+1$, for each $n\in E^0$.  Then 
for $n\in\z$ we may let $N=[n,\infty)$ and $F_k=\emptyset$ for $k\in 
N$.  These are the only elements of $\i(E)$.
\smallskip
We again specialize to directed trees.
\smallskip\noindent
{\bf Definition 1.25.}\enspace  Let $T$ be a directed tree.  Let
$\j\equiv\j(T)$ be the set of open subsets of 
$\partial T$.  We let $\j$ be ordered by inclusion.
We define maps $U\subdel:\i\to\j$ and $L:\j\to\i$ as 
follows.
$$\eqalign{U\subdel(N,F)&= \bigcup_{u\in N}
\bigl[V(u;F_u)\bigr]\subdel\cr\cr
L(W)&= (N,F)\equiv\bigl(N(W),F(W)\bigr),\cr}$$
where
$$\displaylines{\quad
N=\bigl\{u\in T^0\setminus\Sigma(T)\bigm| \exists 
F\subseteq\Delta_1(u)\hbox{\ finite\ with\ }\bigl[V(u;F)\bigr]\subdel 
\subseteq W\bigr\}\hfill\cr
\hskip1true in \cup\  \bigl\{u\in T^0\bigm| 
\bigl[V(u)\bigr]\subdel
\subseteq W\bigr\},\hfill\cr\cr
\quad F_u= \bigcap\bigl\{F\subseteq\Delta_1(u)\bigm|\bigl[V(u;F)
\bigr]\subdel\subseteq W\bigr\},\quad u\in N.\hfill\cr}$$
\theorem Lemma 1.26. The above definitions make sense, 
in that $U\subdel(N,F)\in\j$ 
and $L(W)\in\i$ whenever $(N,F)\in\i$ and $W\in\j$.\endtheorem
\pf Left to the reader. \qed
\smallskip\noindent
{\bf Remark 1.27.}\enspace  It is immediate that the 
maps $L$ and $U\subdel$ are 
order-preserving.
\theorem Theorem 1.28. The maps $L$ and $U\subdel$ are inverses.\endtheorem
\pf We first show that $U\subdel\circ L$ is the identity on $\j$.  Let 
$W\in\j$, and let $L(W)=(N,F)$.  Let $u\in N$.
Then $\bigl[V(u;F_u)\bigr]\subdel$ is one of the 
defining sets of $U\subdel\bigl(L(W)\bigr)$.  By the definition of $L(W)$ we 
have $\bigl[V(u;F_u)\bigr]\subdel\subseteq W$.  
It follows that $U\subdel\bigl(L(W)\bigr)\subseteq W$.

For the reverse inclusion, note that since $W$ is open, it is a union 
of basic open sets.  Let $\bigl[V(u;F)\bigr]\subdel\subseteq W$ 
be one, where $F$ is a finite subset of $\Delta_1(u)$.
We may assume that $F$ is the smallest subset of 
$\Delta_1(u)$ for which the containment holds.  Now, if 
$u\in\partial T$, then $u\in N$.  It follows that 
$\bigl[V(u;F)\bigr]\subdel\subseteq U\subdel\bigl(L(W)\bigr)$.  The 
same argument holds if $u\notin\partial T$ but $u\in N$.  
Finally, we suppose that $u\notin N$.  Then we must have 
$u\not\in\partial T$ and
$\bigl[V(u)\bigr]\subdel\not\subseteq W$.  Then $\Delta_1(u)$ is 
finite.  We have
$$\eqalign{W&\supseteq \big[V(u;F)\bigr]\subdel\cr
&=\bigcup_{e\in\Delta_1(u)\setminus 
F}\bigl[V\bigl(t(e)\bigr)\bigr]\subdel.\cr}$$
Since $\bigl[V\bigl(t(e)\bigr)\bigr]\subdel\subseteq W$ for such $e$, 
we have $t(e)\in N$.  Then 
$\bigl[V\bigl(t(e)\bigr)\bigr]\subdel\subseteq U\subdel\bigl(L(W)\bigr)$.  It 
now follows that $\big[V(u;F)\bigr]\subdel\subseteq  
U\subdel\bigl(L(W)\bigr)$.  Therefore $W\subseteq U\subdel\bigl(L(W)\bigr)$.

We next show that $L\circ U\subdel$ is the identity on $\i$.  Let 
$(N,F)\in\i$, and let $L\bigl(U\subdel(N,F)\bigr)=(M,G)$.  Let $u\in N$.  
Then $\bigl[V(u;F_u)\bigr]\subdel\subseteq U\subdel(N,F)$, (where 
$F_u=\emptyset$ if $\Delta_1(u)$ is finite).  Hence $u\in M$, and so 
$N\subseteq M$.  It also 
follows that $G_u\subseteq F_u$.

We now prove $N\supseteq M$.  Let $w_0\in T^0$ and $w_0\notin N$.  We 
distinguish two cases.
\smallskip\noindent
{\it Case (i).\/} There exists $p_0\in T^*_+$ with $w_0=o(p_0)$, 
$t(p_0)\in\partial T$, and $t(p_0)\notin N$.
\smallskip\noindent
In this case, we claim first that $\bigl[t(p_0)\bigr]\notin U\subdel(N,F)$.  
For suppose 
otherwise.  Then there is $u\in N$ with
$$\bigl[t(p_0)\bigr]\in\bigl[V(u;F_u)\bigr]\subdel.$$
But then $t(p_0)\in V(u;F_u)\subseteq N$, a contradiction.  
This establishes the claim.

Now, if $w_0\in\partial T$,
 then letting $p_0=w_0$ above we obtain $w_0\notin U\subdel(N,F)$.  
Hence $w_0\notin M$.  On the other hand, if $\Delta_1(w_0)$ 
is finite and nonempty, 
the above shows that $t(p_0)\notin U\subdel(N,F)$, while $t(p_0)\in\bigl[ 
V(w_0)\bigr]\subdel$.  Hence  $\bigl[V(w_0)\bigr]\subdel\not\subseteq 
U\subdel(N,F)$, so $w_0\notin M$.
\smallskip\noindent
{\it Case (ii).\/} For all $p_0\in T^*_+\cap o\inv(w_0)$, if 
$t(p_0)\in\partial T$ then $t(p_0)\in N$.
\smallskip\noindent
In particular, since $w_0\notin N$, we have that $\Delta_1(w_0)$ is 
finite and nonempty.  Applying Definition 1.23 (iii) to $w_0$, we have two 
possibilities:
\smallskip
\item{(iia)} There is $w_1\in V(w_0)\setminus\{w_0\}$ with $w_1\in N$ 
and $F_{w_1}\ne\emptyset$, (and hence $w_1\in\partial T$).
\smallskip
\item{(iib)} For every $w\in V(w_0)\setminus\{w_0\}$ we have $w\in N 
\implies F_w=\emptyset$, and there exists $w_1\in 
V(w_0)\setminus\{w_0\}\setminus N$.
\smallskip\noindent
We will treat case (iib) first.  By the assumption of case (ii), we 
must have $\Delta_1(w_1)$ finite. If $w_1$ leads to a vertex in 
$\partial T\setminus N$, then so does $w_0$, contradicting the 
assumption of case (ii).  By the hypothesis of case (iib), $w_1$ must 
also fall under case (iib). So we obtain $w_2\in 
V(w_1)\setminus\{w_1\} \setminus N$.  Inductively, we obtain a 
sequence $w_0$, $w_1$, $\ldots$ in $T^0\setminus N$ with 
$\Delta_1(w_j)$ finite for all $j$, and a path $q\in T^\infty_+\cap 
o\inv(w_0)$ passing through all of the $w_j$.

In case (iia), we have $w_1\in V(w_0)\setminus\{w_0\}$ with $w_1\in 
N$ and $F_{w_1}\ne\emptyset$.  Let $e_1\in F_{w_1}$. If $t(e_1)\in N$, 
then Definition 1.23 (ii) (b) implies that $F_{t(e_1)}\ne\emptyset$.  
Hence if $e\in F_{t(e_1)}$, then $e\in F_{o(e)}$.  If $t(e_1)\notin 
N$, then since we are in case (ii) (for $w_0$), 
$\Delta_1\bigl(t(e_1)\bigr)$ is finite.  We then have case (iia) or (iib) 
for $t(e_1)$.  If case (iib) applies, then as before, we may extend 
from $t(e_1)$ an infinite directed path having infinitely many 
vertices not in $N$.  If case (iia) applies, then we may extend from 
$t(e_1)$ a finite directed path whose last edge satisfies $e\in 
F_{o(e)}$, and we repeat this process.  

Thus in all cases,
either we have extended from $w_0$ an infinite directed path 
having infinitely many vertices not in $N$, or we have extended from 
$w_0$ an infinite directed path having infinitely many edges 
satisfying $e\in F_{o(e)}$.

Let $q$ be this path. We claim that $[q]\notin U\subdel(N,F)$.  For if it were, 
there would be $u\in N$ such that 
$[q]\in\bigl[V(u;F_u)\bigr]\subdel$.  It then follows that
there is a path $p\in T^*\cap o\inv(u)$ such that $o(p)=w_0$ 
and $q=pf_1f_2\cdots$ in $T^\infty$, where $f_j\in T^1_+$ (there may 
be some cancellations in the product $pf_1f_2\cdots$).  Since $f_1\in 
\Delta_1(u)\setminus F_u$, we have $V\bigl(t(f_1)\bigr)\subseteq N$, 
and $F_v=\emptyset$ for all $v\in V\bigl(t(f_1)\bigr)$, by Definition 
1.23 (ii) (a).  Since at most finitely many $f_j$ can have cancelled, we 
have contradicted the construction of $q$.  Therefore 
$[q]\notin U\subdel(N,F)$, and hence 
$\bigl[V(w_0)\bigr]\subdel\not\subseteq U\subdel(N,F)$.  Therefore $w_0\notin 
M$, since $\Delta_1(w_0)$ is finite and nonempty,
and we have shown that $N=M$. 

Finally, we show that $G=F$.  We already know that $G_u\subseteq F_u$ 
for all $u\in N$.  Suppose that there is $u\in N$ such that $G_u\ne 
F_u$.  Let $e_1\in F_u\setminus G_u$.  Then $t(e_1)\in N$,  
$G_{t(e_1)}=\emptyset$, and $F_{t(e_1)}\ne\emptyset$ (by Definition 1.23 
(iii)).  We choose $e_2\in F_{t(e_1)}$.  Since $e_2\notin G_{t(e_1)}$,
$t(e_2)\in N$,  
$G_{t(e_2)}=\emptyset$, and $F_{t(e_2)}\ne\emptyset$.  Inductively, we 
obtain $q=e_1e_2\cdots\in T^\infty_+$ such that
$$\eqalign{o(q)&=u\cr G_{t(e_j)}&=\emptyset,\hbox{\ for\ all\ }j\cr
e_j&\in F_{o(e_j)},\hbox{\ for\ all\ }j.\cr}$$
It is clear that $[q]\in U\subdel(N,G)$.  In the same way as in the previous
proof, it follows that $[q]\notin U\subdel(N,F)$.  But since $(N,G) = 
U\subdel\Bigl(L\bigl(U\subdel(N,F)\bigr)\Bigr)$, it follows from the fact that 
$U\subdel\circ L=\hbox{id}\;_\j$ that $U\subdel(N,G)=U\subdel(N,F)$.  
This contradiction 
finishes the proof. \qed

We next show that a closed subset of the boundary of a directed tree 
corresponds to an algebra as described in Definition 1.16, corresponding 
to a certain subforest of the tree.  
Again, since we will need to apply this idea to graphs later, we 
define it now  for directed graphs.
\smallskip\noindent
{\bf Definition 1.29.} \enspace Let $E$ be a directed graph, and let 
$(N,F)\in\i(E)$ be an invariant of $E$.
\item{(i)} $R(N,F)=\bigl\{u\in N\bigm|F_u\ne\emptyset\bigr\}$.
\smallskip
\item{(ii)} $E(N,F)$ is the subgraph of $E$ given by
\itemitem{} $E(N,F)^0=(E^0\setminus N)\cup R(N,F)$
\itemitem{} $E(N,F)^1_+=E^1_+\cap o\inv\bigl(E(N,F)^0\bigr) \cap t\inv 
\bigl(E(N,F)^0\bigr)$.
\smallskip
\item{(iii)} $S(N,F) = R(N,F) \cup \bigl\{u\in E^0 \setminus N\bigm| 
0< \#\Delta_1(u)<\infty\bigr\}$.

\theorem Proposition 1.30. Let $T$ be a directed tree and let 
$(N,F)\in\i(T)$.  Then
$$ C_0\bigl(\partial 
T\bigr)/C_0\bigl(U\subdel(N,F)\bigr)\cong
A\bigl(T(N,F),S(N,F)\bigr).$$
\endtheorem

\pf We define open subsets of $A(T)\hat$ by
$$\eqalign{U(N,F)&= \bigcup_{u\in N}\bigl[V(u;F_u)\bigr]\cr
P(N,F)&=\bigcup_{u\in N\setminus R(N,F)}\bigl[V(u)\bigr].\cr}$$
Thus $U\subdel(N,F)=U(N,F)\cap\partial T$ and $U(N,F)=P(N,F) \cup 
R(N,F)$.

Next we note that 
$$A\bigl(T(N,F)\bigr)\cong A(T)\bigm/ C_0\bigl(P(N,F)\bigr).$$
For the map
$$V_T(u)\mapsto\cases{V_{T(N,F)}(u),&if $u\in T(N,F)^0$\cr
\emptyset,&if $u\in N\setminus R(N,F)$\cr}$$
extends to a Boolean ring homomorphism by Lemma 1.5.  Then the kernel 
of the corresponding $*$-homo\-mor\-phism of $A(T)$ onto 
$A\bigl(T(N,F)\bigr)$ is generated by the characteristic functions 
of $\bigl\{\bigl[V_T(u)\bigr]\bigm|u\in N\setminus R(N,F)\bigr\}$.  
Hence the kernel equals $C_0\bigl(P(N,F)\bigr)$.

Now,
$$\eqalign{{C_0(\partial T)\over C_0\bigl(U\subdel(N,F)\bigr)}&=
{A(T)\bigm/ C_0\bigl(\Sigma(T)\bigr)\over C_0\bigl(U(N,F)\cap\partial 
T\bigr)}\cr
&= {A(T)\bigm/C_0\bigl(\Sigma(T)\bigr)\over C_0\bigl(U(N,F)\bigr) 
\bigm/C_0\bigl(U(N,F)\cap\Sigma(T)\bigr)}\cr
&\cong {A(T)\over C_0\bigl(U(N,F)\cup\Sigma(T)\bigr)}\cr
&\cong {A(T)\bigm/C_0\bigl(P(N,F)\bigr)\over C_0\bigl(U(N,F)\cup 
\Sigma(T)\bigr)\bigm/ C_0\bigl(P(N,F)\bigr)}.\cr}$$

\noindent
We have already noted that the numerator is $A\bigl(T(N,F)\bigr)$.  
As for the denominator, we have
$$\eqalign{C_0\bigl(U(N,F)\cup\Sigma(T)\bigr)\bigm/C_0\bigl( 
P(N,F)\bigr) &= C_0\bigl(R(N,F)\cup\Sigma(T)\cup P(N,F)\bigr) \bigm/ 
C_0\bigl(P(N,F)\bigr)\cr
&\cong C_0\Bigl(\bigl(R(N,F)\cup\Sigma(T)\bigr)\setminus P(N,F)\Bigr)\cr
&=C_0\Bigl(R(N,F)\cup\bigl(\Sigma(T)\setminus N\bigr)\Bigr)\cr
&= C_0\Bigl(R(N,F)\cup\bigl(\Sigma(T)\cap T(N,F)^0\bigr)\Bigr)\cr
&= C_0\bigl(S(N,F)\bigr).\cr}$$
Hence $C_0(\partial T)\bigm/ C_0\bigl(U\subdel(N,F)\bigr)\cong 
A\bigl(T(N,F),S(N,F)\bigr)$. \qed

\bigskip\noindent
\centerline{\bf 2.  Groupoids and $C^*$-algebras for Directed Graphs.}
\smallskip
We now let $E$ be a directed graph.
  The set $E^*$ of (undirected) finite paths has two 
structures that we will use.

\smallskip\noindent
{\bf Definition 2.1.}\enspace  We let $G=G(E)$ denote the set $E^*$ with the 
following groupoid structure:  
$$\displaylines{
\qquad G^0=E^0\hfill\cr
\qquad r(p)=o(p),\ s(p)=t(p), \hbox{\ and\ }
(e_1\cdots e_n)^{-1}=\overline{e_n}\cdots 
\overline{e_1},\hbox{\ for\ }p=e_1\cdots e_n
\in G\hfill\cr 
\qquad G^{(2)}=\bigl\{(p,q)\in G\times 
G\bigm|t(p)=o(q)\bigr\}\hfill\cr
\qquad \hbox{multiplication\ is\ given\ by\ 
concatenation,\ followed\ by}\hfill\cr
\qquad\qquad\hbox{the\ removal\ of\ pairs\ of\ the\ 
form\ } 
(e,\overline e)\hfill\cr}$$

\smallskip\noindent
{\bf Definition 2.2.}\enspace   We let $\tildee$ denote the following graph:
$$\displaylines{
\qquad\tildee^0=E^*\hfill\cr
\qquad\tildee^1=\bigl\{(p,q)\in E^*\times E^*\bigm|(p^{-1},q)\in 
G^{(2)}\hbox{\ and\ }p^{-1}q\in E^1\bigr\}\hfill\cr
\qquad\tildee^1_+=\bigl\{(p,q)\in\tildee^1\bigm| p^{-1}q\in 
E^1_+\bigr\}\hfill\cr
\qquad o(p,q)=p,\ t(p,q)=q,\ \overline{(p,q)}=(q,p)\hfill\cr}$$

\smallskip\noindent
{\bf Definition 2.3.}\enspace  The map $\lambda:\tildee^{**}
\to E^0$ is defined by
$$\eqalign{
\lambda(p)&=o(p),\quad p\in \tildee^0=E^*,\cr
\lambda(\mu)&=\lambda(o(\mu)),\quad\mu\in \tildee^n
\hbox{\ for\ }1\le 
n\le \infty.\cr}$$

\theorem Lemma 2.4.  With $\lambda$ as fiber map,
 $\tildee$ is a bundle of directed trees.\endtheorem

\pf It is easy to check that $p$, $q\in\tildee^0$
 are connected by a 
path if and only if $\lambda(p)=\lambda(q)$, and 
that in this 
case the path 
connecting them is unique.\qed

\smallskip
For $x\in\tildee^\infty_+/\sim$ we let $\lambda(x)$ 
denote the image 
under $\lambda$ of any representative of $x$.
The groupoid $G$ acts on $\tildee$ in the obvious way:  
if $\alpha\in 
G$ and $p\in\tildee^0$ satisfy $s(\alpha)=\lambda(p)$ 
then $\alpha\cdot p$ is obtained by multiplication in $G$.  
If $\alpha\in 
G$ and $(p,q)\in\tildee^1$ satisfy $s(\alpha)=\lambda(p,q)$ 
then $\alpha\cdot(p,q)=(\alpha\cdot p,\alpha\cdot q)$.  
It is clear 
that the action of $G$ on $\tildee$ preserves 
the direction.

\smallskip\noindent
{\bf Definition 2.5.}\enspace Let $E$ be a directed graph.  We define 
$i:I(E)\to I(\tildee)$ as follows.  For $(N,F)\in I(E)$, let 
$i(N,F)=(\tilden,\tildef)$, where
$$\eqalign{\tilden&=t\inv(N)\cr
\tildef_p&=\bigl\{(p,pe)\bigm|e\in F_{t(p)}\bigr\}.\cr}$$

\smallskip\noindent
{\bf Remark 2.6.}\enspace The map $i$ is clearly injective and 
order-preserving.

\theorem Lemma 2.7. Let $E$ be a directed graph, and
let $(\Lambda,\Phi)\in I(\tildee)$.  The following are equivalent:
\item{(i)} $(\Lambda,\Phi)\in i\bigl(I(E)\bigr)$.
\item{(ii)} $U\subdel(\Lambda,\Phi)$ is $G(E)$-invariant.
\item{(iii)} $\Lambda = t\inv\bigl(t(\Lambda)\bigr)$, and for 
$p\in\Lambda$, $\Phi_p=p\Phi_{t(p)}$.
\endtheorem

\pf $\bigl(\hbox{(i)}\Rightarrow\hbox{(ii)}\bigr):\quad$ 
Let $(\Lambda,\Phi)=(\tilden,
\tildef)$ 
for some $(N,F)\in I(E)$.  Note that 
$$\tildef_p = p\bigl\{\bigl(t(p),e\bigr)\bigm|e\in F_{t(p)}\bigr\},$$
(using the action of $G(E)$ on $\tildee$), and hence that
$$\bigl[V_\tildee(p;\tildef_p)\bigr]\subdel = 
p\bigl[V_\tildee(t(p);\tildef_{t(p)})\bigr]\subdel.$$
Since $\tilden$ is $G(E)$-invariant, it follows that 
$U\subdel(\tilden,\tildef)$ is $G(E)$-invariant.
\smallskip\noindent
$\bigl(\hbox{(ii)}\Rightarrow\hbox{(iii)}\bigr):\quad$ 
Let $p\in\Lambda$.  Then
$$\bigl[V_\tildee(p;\Phi_p)\bigr]\subdel\subseteq U\subdel 
(\Lambda,\Phi).$$
If $t(q)=t(p)$, then by $G(E)$-invariance we have
$$\bigl[V_\tildee(q;qp\inv\Phi_p)\bigr]\subdel\subseteq 
U\subdel(\Lambda,\Phi).$$
By Theorem 1.28 we have
\item{(1)} $q\in\Lambda$,
\item{(2)} $qp\inv\Phi_p\supseteq\Phi_q$.
\par\noindent
It follows from (1) that $\Lambda=t\inv\bigl(t(\Lambda)\bigr)$.  It 
follows from (2) and symmetry that $p\inv\Phi_p=q\inv\Phi_q$.  In 
particular, $p\inv\Phi_p=\Phi_{t(p)}$.
\smallskip\noindent
$\bigl(\hbox{(iii)}\Rightarrow\hbox{(i)}\bigr):\quad$ 
Let $N=t(\Lambda)$.  Then 
$\Lambda=t\inv\bigl(t(\Lambda)\bigr)=t\inv(N)=\tilden$.  Since 
vertices of $E$ are paths of length zero, if $p\in\Lambda$ then 
$t\bigl(t(p)\bigr)=t(p)\in\Lambda$.  Thus 
$t(p)\in t\inv\bigl(t(\Lambda)\bigr)=\Lambda$.  Hence 
$N\subseteq\Lambda$. For $u\in N$ let 
$$F_u=\bigl\{e\in\Delta_{1,E}(u)\bigm|(u,e)\in\Phi_u\bigr\}.$$
The conditions that $(N,F)\in I(E)$ follow from those that 
$(\Lambda,\Phi)\in I(\tildee)$.
By hypothesis, then, if $t(p)=u$ we have
$$\eqalign{\Phi_p&=p\Phi_u\cr
&=\bigl\{(p,pe)\bigm|e\in F_u\bigr\}\cr
&=\tildef_p.\cr}$$
Thus $(\Lambda,\Phi)=(\tilden,\tildef)=i(N,F)$. \qed

\theorem Theorem 2.8.  Let $E$ be a directed graph.  The map 
$U\subdel\circ i$ is an order-preserving one-to-one correspondence 
between $I(E)$ and the collection of open $G(E)$-invariant subsets 
of $\partial\tildee$.\endtheorem

\pf This follows from Theorem 1.28 and Lemma 2.7. \qed

\smallskip
We now prepare for the definition of the groupoids of a directed graph.

\smallskip\noindent
{\bf Definition 2.9.}\enspace  Let $E$ be a directed graph, and let $x\in 
A(\tildee)\hat$.  We 
define a function $\underline\cdot:A(\tildee)\hat\to E^{**}$
such that $o(\underline{x})=\lambda(x)$ as 
follows.  If $x\in\tildee^0$, then $x\in E^*$ and we 
let $\underline 
x=x$.  If $x\in \tildee^\infty_+/\sim$, then there is 
a unique $\omega\in 
\tildee^\infty$ with $o(\omega)=\lambda(x)$ and such 
that some tail of 
$\omega$ is in (the equivalence class) $x$. Then we 
may write $\omega$ 
as
$$\omega=(p_0,p_1)(p_1,p_2)\cdots,$$
where $p_0=\lambda(x)$ and for each $i$, $\ell(p_i)=i$.
We define $\underline x=(p_0^{-1}p_1)(p_1^{-1}p_2)
\cdots\in E^\infty$.

\theorem Lemma 2.10.  The map $$\underline\cdot:
A(\tildee)\hat\to 
E^{**}$$ is one-to-one and has range 
$\bigl\{p\in E^{**}\bigm| \hbox{\ 
some\ tail\ of\ $p$\ is\ directed}\bigr\}$.  
For $(\alpha,x)\in 
G*A(\tildee)\hat$, $\underline{\alpha x}=\alpha\underline{x}$.
\endtheorem

\pf The map is a bijection from $\tildee^0\to E^*$.  For 
$x\in\tildee^\infty_+/\sim$, if 
$\underline x=f_1f_2\cdots$, then 
we may construct an infinite path in the equivalence 
class $x$ by
letting $p_i=f_1f_2\cdots f_i$ (and $p_0=o(f_1)$), and setting 
$\omega=(p_0,p_1)(p_1,p_2)\cdots$.  This defines a left 
inverse for the map $\underline\cdot$, proving 
it is one-to-one.  The 
description of the range is clearly true.  
To check the equivariance 
of the $G$-action, let $(\alpha,x)\in G*A(\tildee)\hat$.
It is clear that equivariance holds if $x\in\tildee^0$.  So let 
$x\in\tildee^\infty_+/\sim$, and let
$\omega=(p_0,p_1)(p_1,p_2)\cdots$ be as in Definition 2.9.  
Let $\alpha=e_1\cdots e_n$, and suppose that the product 
$\alpha\underline x$ involves $k$ cancellations.  Then
$$\eqalign{
\alpha \underline x&=(e_1\cdots 
e_n)\cdot\bigl(p_1(p_1^{-1}p_2)\cdots\bigr)\cr
&=e_1\cdots e_{n-k}(p_k^{-1}p_{k+1})\cdots.\cr}$$
But the unique element of $\tildee^\infty$ with origin equal to 
$\lambda(\alpha x)$ and with a tail representing $\alpha x$ is
$$\Bigl(\bigl(o(e_1),e_1\bigr)\bigl(e_1,e_1e_2\bigr)
\cdots\bigl(e_1\cdots e_{n-k},
e_1\cdots e_{n-k}(p_k^{-1}p_{k+1})\bigr)\cdots\Bigr)$$
Then it is clear that $\underline{\alpha x}=
\alpha\underline{x}$. \qed

\smallskip\noindent
{\bf Definition 2.11.}\enspace  
For $S\subseteq\Sigma(E)$ we let $\tildes=t^{-1}(S) = \{p\in 
\tildee^0\mid t(p)\in S\}$.


\smallskip\noindent
{\bf Remark 2.12.}\enspace   If $S\subseteq\Sigma(E)$, then $\tildes$ 
is open and $G$-invariant.

\smallskip\noindent
Let $S\subseteq\Sigma(E)$.  We define
$$\eqalign{
A(\tildee,\tildes)&=\bigoplus_{v\in E^0}A\bigl(\lambda^{-1}(v), 
\tildes\cap\lambda^{-1}(v)\bigr)\cr
\partial(\tildee,\tildes)&=A(\tildee,\tildes)\hat=
\bigcup_{v\in E^0} \partial\bigl(
\lambda^{-1}(v), 
\tildes\cap\lambda^{-1}(v)\bigr)\cr}$$
where the union is given the inductive topology.  

\smallskip\noindent
{\bf Definition 2.13.}\enspace  Let $E$ be a directed graph.  Let 
$S\subseteq\Sigma(E)$.
The transformation groupoid 
$G*\partial(\tildee,\tildes)$ 
(with the relative product topology) is called an {\it extended 
Toeplitz graph
groupoid} of $E$ (determined by the choice of $S$). 

\smallskip
We recall the groupoid structure of $G*\partial(\tildee,\tildes)$:
$$\eqalign{\bigl(G*\partial(\tildee,\tildes)\bigr)^0&=
\partial(\tildee,\tildes)\cr
s(\alpha,x)&=x\cr
r(\alpha,x)&=\alpha\cdot x\cr
(\alpha,\beta x)(\beta,x)&=(\alpha\beta,x)\cr
(\alpha,x)^{-1}&=(\alpha^{-1},\alpha x).\cr}$$
Since $G$ is discrete it is clear that $s$ and $r$ are local 
homeomorphisms, and hence that $G*\partial(\tildee,\tildes)$
is $r$-discrete.
We call the \cstar-algebra of $G*\partial(\tildee,\tildes)$ 
an {\it extended Toeplitz graph algebra}.

\theorem Lemma 2.14.  Let $S\subseteq\Sigma(E)$, and let 
$D=G*\partial(\tildee,\tildes)$.  The nondegenerate 
$*$-represen\-ta\-tions of $C_c(D)$ are in one-to-one 
correspondence with 
the pairs $(\pi,U)$, where $\pi$ is a
 nondegenerate representation of 
$C_c(\partial\bigl(\tildee,\tildes)\bigr)$ on a Hilbert space $H$, 
and $U:G\to L(H)$ is a representation of $G$ 
by partial isometries such that 
$$\eqalign{U_\alpha^*U_\alpha
&=\pi\bigl(\Chi_{\lambda^{-1}(s(\alpha))}\bigr),\cr 
U_\alpha U_\alpha^*&=
\pi\bigl(\Chi_{\lambda^{-1}(r(\alpha))}\bigr),\cr
U_\alpha\pi\bigl(\Chi_{V(p)}\bigr)U_\alpha^*
&=\pi\bigl(\Chi_{V(\alpha p)}\bigr),\cr}$$
where $V(\alpha p)=\emptyset$ if $s(\alpha)\not=\lambda(p)$, and
$\pi$ also denotes the extension of $\pi$ to the 
multiplier algebra 
of $C_0\bigl(\partial(\tildee,\tildes)\bigr)$.\endtheorem

\pf If $\sigma:C_c(D)\to L(H)$ is a nondegenerate $*$-representation, 
then $\sigma\restrictedto{C_c(D^0)}$ extends to a nondegenerate 
$*$-representation $\pi:C_0\bigl(\partial(\tildee,\tildes)\bigr)\to 
L(H)$.  Since $\lambda^{-1}(u)$ is a clopen subset of 
$\partial(\tildee,\tildes)$ for $u\in E^0$, the extension of $\pi$ 
to $C_b(\partial\bigl(\tildee,\tildes)\bigr)$ can be applied to 
$\Chi_{\lambda^{-1}(u)}$.
\smallskip
For $\alpha\in G$,
$$\bigl\{\sigma(\Chi_{\{\alpha\}\times V})\bigm| V\subseteq 
\lambda^{-1}\bigl(s(\alpha)\bigr) \hbox{\ is\ compact-open} \bigr\}$$ 
is a coherent net of partial isometries, hence converges in
the strong operator topology to a partial isometry $U_\alpha\in 
L(H)$.  It is clear that $U_\alpha^*U_\alpha$ is orthogonal to 
$\pi(\Chi_{\lambda^{-1}(u)})$ if $u\not=s(\alpha)$.  
By the nondegeneracy of $\pi$, it follows that $U_\alpha^*U_\alpha 
=\pi(\Chi_{\lambda^{-1}\bigl(s(\alpha)\bigr)})$.  Similarly,
$U_\alpha U_\alpha ^*
=\pi(\Chi_{\lambda^{-1}\bigl(r(\alpha)\bigr)})$.  The remaining 
properties of the $\{U_\alpha\}$ follow from the definition.
\smallskip
Conversely, if $(\pi,U)$ are given, we define $\sigma:C_c(D)\to L(H)$ 
by
$$\sigma\bigl(\Chi_{\{\alpha\}\times\lambda^{-1}
\bigl(s(\alpha)\bigr)}\cdot f\bigr)= U_\alpha\pi(f)$$
whenever $f\in C_c(D^0)$ and 
$\supp(f)\subseteq\lambda^{-1}\bigl(s(\alpha)\bigr)$.  Then $\sigma$ 
is a linear map, and it is easy to verify that $\sigma$ is 
multiplicative and $*$-preserving by checking on elements of the 
form $\Chi_{\{\alpha\}\times V}$ for $V\subseteq\lambda^{-1} 
\bigl(s(\alpha)\bigr)$ compact-open. \qed

\theorem Corollary 2.15.  Every $*$-representation of 
$C_c\bigl(G*\partial(\tildee,\tildes)\bigr)$ is bounded. \qed
\endtheorem
\noindent
(This is true for any $r$-discrete groupoid --- 
see \bib\quiggsieben, 
Proposition 3.2.)

\smallskip
To obtain more familiar objects we will define natural 
transversals 
in the extended Toeplitz graph groupoids.
Let $S\subseteq\Sigma(E)$.  We define
$X(E,S)\subseteq\partial(\tildee,\tildes)$ by
$$X(E,S)=\bigl\{x\in \partial(\tildee,\tildes)\bigm| \underline 
x\in E^{**}_+\bigr\}.$$
(Thus $X(E,S)$ is the set of points in the spectrum of 
$A(\tildee,\tildes)$ that are represented by directed paths in 
$E$.)  In the special case $S=\Sigma(E)$ we let 
$X(E)=X(E,\Sigma(E))$, the set of all points of $\partial\tildee$ 
represented by directed paths in $E$.

\theorem Lemma 2.16.  $X(E,S)$ is a clopen transversal in
$G*\partial(\tildee,\tildes)$.\endtheorem

\pf Let $x\in X(E,S)$.  Then $\lambda(x)\in\tildee^0$, so we may 
consider the open set $[V_\tildee(\lambda(x))]$ in $A(\tildee)\hat$.  
Since $x\in X(E,S)$, $x\in [V_\tildee(\lambda(x))]$.  Moreover, 
for any $y\in 
[V_\tildee(\lambda(x))]$, $y$ consists entirely of directed edges 
of $\tildee$, 
and hence $\underline y\in E^{**}_+$.  Thus $[V_\tildee(\lambda(x))]
\cap\partial(\tildee,\tildes)\subseteq 
X(E,S)$, and so $X(E,S)$ is open.
\smallskip
Now suppose that $x\in\partial(\tildee,\tildes)\setminus
 X(E,S)$.  Letting $\underline 
x=e_1e_2\cdots $, then $e_i\not\in E^1_+$ for some $i$.  Let 
$p=e_1\cdots e_i\in\tildee^0$.  Then $x\in [V(p)]$ and $[V(p)]\cap 
X(E,S)=\emptyset$.  Thus the complement of $X(E,S)$ is open.
\smallskip
To see that $X(E,S)$ is a transversal, let $x\in 
\partial(\tildee,\tildes)$.  Let $\underline 
x=e_1e_2\cdots $.  Then there is $n$ with $e_i\in E^1_+$ for $i\ge 
n$.  Let $\alpha=e_1\cdots e_{n-1}$.  Then $(\alpha^{-1},x)$ has 
source $x$ and range in $X(E,S)$.\qed

\smallskip\noindent
{\bf Definition 2.17.}\enspace   Let $E$ be a directed graph.
The {\it Toeplitz graph groupoids} are the 
restrictions of the extended Toeplitz graph groupoids to the 
transversals $X(E,S)$, and the {\it Toeplitz graph algebras} are 
their \cstar-algebras:

$$\eqalign{
\t\g(E,S)&=G(E)*\partial(\tildee,\tildes)\restrictedto{X(E,S)}\cr
\t\oh(E,S)&=C^*\bigl(\t\g(E,S)\bigr),\cr}$$
for $S\subseteq\Sigma(E)$.  (We let $\t\oh(E)\equiv\t\oh(E,\emptyset)$.)

\smallskip\noindent
The {\it graph groupoid} and 
{\it graph algebra}
of $E$ are special cases obtained when  $S=\Sigma(E)$:
$$\eqalign{
\g(E)&=\g(E,X(E))=G(E)*\partial
\tildee\restrictedto{X(E)}\cr
\oh(E)&=C^*\bigl(\g(E)\bigr).\cr}$$

\noindent
{\bf Remark 2.18.}\enspace  In the case of a graph $E$ with no sinks, the graph 
groupoid defined above is exactly the usual notion of infinite 
directed paths ``modulo shift-tail equivalence with lag''.  In the 
case of a graph that is not row-finite, one must include the finite 
directed paths ending at vertices with infinite exit valence (as we 
have done).  This idea was used in an ad-hoc way by previous authors 
(\bib{\renault,\kprr,\szymanskizhang}).

\smallskip
Next we will characterize representations of $\t\oh(E,S)$ by means of 
generators and relations (cf.  \bib{\kpr,\flr}).

\theorem Theorem 2.19. The nondegenerate representations of $\t\oh(E)$ 
are in one-to-one correspondence with the families of operators 
$\{S_e\}_{e\in E^1_+}$ and $\{P_u\}_{u\in E^0}$ such that
\smallskip
\item{(i)} The $\{S_e\}$ are partial isometries and the $\{P_u\}$ are 
projections,
\item{(ii)} $u\not=v\implies P_uP_v=0$,
\item{(iii)} $\sum P_u=1$ in the strong operator topology,
\item{(iv)} $S_e^*S_e=P_{t(e)}$,
\item{(v)} $P_u\ge\sum\{S_fS_f^*\mid o(f)=u\}$.
\endtheorem

\pf If $\pi:\t\oh(E)\to L(H)$ is nondegenerate, let
$$\eqalign{S_e
&=\pi\bigl(\Chi_{\{e\}\times V_\tildee(t(e))}\bigr)\cr 
P_u&=\pi\bigl(\Chi_{V_\tildee(u)}\bigr).\cr}$$
It is easy to verify properties (i) - (v).

\smallskip
Conversely, let $\{S_e\}$, $\{P_e\}$
satisfying (i) - (v) be given.  For $p=e_1e_2\cdots 
e_n\in E^*_+$, let $S_p=S_{e_1}S_{e_2}\cdots S_{e_n}$.
We note the following elementary consequences of (i) - (v):
\smallskip
\item{(a)} $e\not=f\implies S_e^*S_f=0$,
\item{(b)} $S_p^*S_p=P_{t(p)}$,
\item{(c)} $S_p^*S_q=0$ if neither of $p$, $q$ extends the other.
\smallskip\noindent
Let $\theta_p=S_pS_p^*$ for $p\in E^*_+$.  We will verify the 
conditions of Lemma 1.5 for (the bundle of trees) $\tildee$.  If 
$(p,q)\in\Delta_{1,\tildee}(p)$, then $q=pf$ for some $f\in E^1_+$.  Then
$$\eqalign{\theta_q&= S_qS_q^*\cr
&=S_pS_fS_f^*S_p^*\cr
&\le S_pS_p^*\cr
&=\theta_p.\cr}$$
Let $(p,q_1)\not=(p,q_2)\in\Delta_{1,\tildee}(p)$.  
Then there are $f_1\not=f_2 \in 
E^1_+$, with $o(f_1)=o(f_2)=t(p)$, such that $q_i=pf_i$.  Then
$$\eqalign{S_{q_1}^*S_{q_2}&=
S_{f_1}^*S_p^*S_pS_{f_2}\cr
&=S_{f_1}^*P_{t(p)}S_{f_2}\cr
&=S_{f_1}^*S_{f_2}\cr
&=0.\cr}$$
Thus by Lemma 1.5, there is a $*$-homomorphism 
$\pi_0:C_0\bigl(X(E,\emptyset)\bigr)\to L(H)$ such that
$$\eqalign{\pi_0\bigl(\Chi_{V(q)}\bigr)
&=\theta_q,\ q\in E^*_+,\cr
\pi_0\bigl(\Chi_{V(u)}\bigr)
&=P_u,\ u\in E^0.\cr}$$
Then $\pi_0$ is nondegenerate, by (iii).  Define 
$\pi:C_c\bigl(\t\g(E)\bigr)\to L(H)$ as follows.  The elements of 
$C_c\bigl(\t\g(E)\bigr)$ can be written uniquely in the form
$$\sum\Chi_{\{\alpha_1\alpha_2^{-1}\}\times V(\alpha_2)}\cdot 
f_{\alpha_1,\alpha_2},$$
where the sum is taken over pairs $(\alpha_1,\alpha_2)\in E^*_+\times 
E^*_+$ for which $\alpha_1\cdot\alpha_2^{-1}$ involves no 
cancellations, and where $f_{\alpha_1,\alpha_2}
\in C_c\bigl(X(E,\emptyset)\bigr)$ is zero except for finitely many
$(\alpha_1,\alpha_2)$.  We define $\pi$ by
$$\pi\bigl(\Chi_{\{\alpha_1\alpha_2^{-1}\}\times V(\alpha_2)}
\cdot f\bigr) = 
S_{\alpha_1}S_{\alpha_2}^*\pi_0(f).$$
This is well-defined, since
$$\eqalign{\Chi_{\{\alpha_1\alpha_2^{-1}\}\times V(\alpha_2)}\cdot f=0
&\implies \supp(f)\subseteq V(\alpha_2)^c\cr
&\implies S_{\alpha_2}^*\pi_0(f)= S_{\alpha_2}^* 
\theta_{\alpha_2}(1-\theta_{\alpha_2})\pi_0(f)=0.\cr}$$
$\pi$ is clearly linear.  It is straight-forward to verify that $\pi$ 
is multiplicative and adjoint-preserving by checking on elements of 
the form $\Chi_{\{\alpha_1\alpha_2^{-1}\}\times V}$ for $V\subseteq 
[V(\alpha_2)]$ compact-open. \qed

\theorem Theorem 2.20.  Let $S\subseteq\Sigma(E)$.
The nondegenerate representations of 
$\t\oh(E,S)$ 
are in one-to-one correspondence with the families of operators 
$\{S_e\}_{e\in E^1_+}$ and $\{P_u\}_{u\in E^0}$ satisfying (i) - (v) 
of Theorem 2.19 such that equality holds in (v) whenever $u\in S$.
\endtheorem

\pf Note that equality in (v) for $u\in S$ means in particular that 
the sum in (v) is finitely nonzero.  Let $\pi$ be a nondegenerate 
representation of $\t\oh(E)$ corresponding to generators and
relations as in Theorem 2.19.    Then
$$\eqalign{
\pi\hbox{\ factors\ through\ }\t\oh(E,S)
&\Longleftrightarrow C_0\bigl(\tildes\cap X(E,S)\bigr) 
\subseteq\ker\pi\cr
&\Longleftrightarrow \bigl(\forall u\in S\bigr) \bigl(\exists 
F\subseteq\Delta_1(u)\hbox{\ finite}\bigr)\bigl( P_u=\sum_{e\in 
F}S_eS_e^*\bigr)\cr
&\Longleftrightarrow \hbox{\ equality\ holds\ in\ (v)\ 
whenever\ }u\in S.\qed\cr}$$

For reference elsewhere (\bib\spielberg), we give a slightly different form of 
the relations defining the algebra $\oh(E)$.  For this purpose we 
will let $D$ denote the set of vertices with infinite exit valence:
$$D=\bigl\{u\in E^0\bigm|\Delta_1(u)\hbox{\ is\ infinite}\bigr\}.$$

\theorem Theorem 2.21.  Let $E$ be a directed graph.  The nondegenerate 
representations of $\oh(E)$ are in one-to-one correspondence with
the families of operators 
$\{S_e : e\in E^1_+\}$ and $\{P_u : u\in E^0\}$ satisfying
\smallskip
\item{(1)} The $\{S_e\}$ are partial isometries and the $\{P_u\}$ are 
projections,
\item{(2)} $u\not=v\implies P_uP_v=0$,
\item{(3)} $S_e^*S_e=P_{t(e)}$,
\item{(4)} $P_{o(e)}S_e=S_e$, if $o(e)\in D$,
\item{(5)} $e\ne f$ and $o(e)=o(f)\in D\implies S_e^*S_f=0$,
\item{(6)} $u\notin D\implies P_u=\sum\{S_fS_f^*\mid o(f)=u\}$,
\item{(7)} $\sum P_u=1$ in the strong operator topology.
\endtheorem

\pf This follows easily from Theorems 2.19 and 2.20. \qed

\smallskip
Next we will define the fundamental cocycle on the (extended) graph 
goupoids.

\theorem Lemma 2.22.  Let $(\alpha,y)\in G*A(\tildee)\hat$.  Then there 
exist unique elements $\beta_1,\beta_2\in G$ and $x\in 
A(\tildee)\hat$ such that
\smallskip\noindent
\halign{\quad\hfil$#$ &#\hfil\cr
(i)&$(\beta_1,x),\ (\beta_2,x)\in G*A(\tildee)\hat$\cr
(ii)&$\alpha=\beta_1\beta_2^{-1}$\cr
(iii)&$y=\beta_2x$\cr
(iv)&The products $\beta_1\cdot \underline x$, $\beta_2\cdot\underline 
x$, $\beta_1\cdot\beta_2^{-1}$ involve no cancellations.\cr}
\endtheorem

\pf For existence, let $\alpha=e_1\cdots e_k$ and $\underline 
y=f_1f_2\cdots$, and let $r\le k$ denote the number of cancellations 
occurring when $\alpha$ and $\underline y$ are multiplied.  We have 
that
$$\eqalign{f_i&=\overline{e_{k-i+1}},\ 1\le i\le r\cr
f_{r+1}&\not=\overline{e_{k-r}}.\cr}$$
Then we may put
$$\eqalign{\beta_1&=e_1\cdots e_{k-r}\cr
\beta_2&=f_1\cdots f_r\cr
\underline x&=f_{r+1}f_{r+2}\cdots.\cr}$$
For uniqueness, let $\beta_1^\prime,\beta_2^\prime$, and $x^\prime$ 
also satisfy conditions (i)-(iv) of the statement of the Lemma.  
Since $\underline y=\beta^\prime\underline x^\prime$, there is $s$ 
such that
$$\eqalign{\beta_2^\prime&=f_1\cdots f_s\cr
\underline x^\prime&=f_{s+1}f_{s+2}\cdots.\cr}$$
If $\beta_1^\prime=e_1^\prime\cdots e_\ell^\prime$, then from 
$\beta_1\cdot\underline x=\beta_1^\prime\cdot\underline x^\prime$
we have
$$e_1\cdots e_{k-r}f_{r+1}f_{r+2}\cdots= e_1^\prime\cdots 
e_\ell^\prime f_{s+1}f_{s+2}\cdots.$$
Suppose $r<s$.  Since paths in a tree are unique (in this case, the 
path from $o(e_1)=o(e_1^\prime)$ to $o(f_{s+1})$), we must have
$$e_1\cdots e_{k-r}f_{r+1}\cdots f_s=e_1^\prime\cdots 
e_\ell^\prime=\beta_1^\prime.$$
But then the product $\beta_1^\prime(\beta_2^\prime)^{-1}$ involves 
the cancellation of $f_{r+1}\cdots f_s$, a contradiction.  The 
assumption $r>s$ leads to a similar contradiction.  Therefore $r=s$, 
and it follows that $\beta_1^\prime=\beta_1$, 
$\beta_2=\beta_2^\prime$, and $x^\prime=x$.\qed

\smallskip\noindent
{\bf Definition 2.23.}\enspace  Let $(\alpha,y)\in G*A(\tildee)\hat$.  The {\it 
standard form} of $(\alpha,y)$ is the triple $(\beta_1,\beta_2,x)$ 
satisfying (i) - (iv) of Lemma 2.22.

\smallskip\noindent
{\bf Definition 2.24.}\enspace  $c:G*A(\tildee)\hat\to\z$ is defined by 
$c(\alpha,y) =\ell(\beta_1)-\ell(\beta_2)$, where 
$(\beta_1,\beta_2,x)$ is the standard form of $(\alpha,y)$.  For 
$S\subseteq\Sigma(E)$, we also let $c$ denote the restriction 
to $G*\partial(\tildee,\tildes)$.

\theorem Lemma 2.25.  $c$ is a continuous homomorphism.\endtheorem

\pf The continuity is clear.  Note that it follows from the 
definition of $c$ that $c$ is zero on units, and that 
$c(\alpha,y)^{-1} =-c(\alpha,y)$.  We will prove that 
$c\bigl((\beta,y)(\alpha,y)^{-1}\bigr)=c(\beta,y)-c(\alpha,y)$.  Let 
$y=g_1g_2\cdots$, $\beta=e_1\cdots e_j 
\overline{g_k}\cdots\overline{g_1}$, $\alpha=f_1\cdots f_p 
\overline{g_q}\cdots\overline{g_1}$, where 
$e_j\not=\overline{g_{k+1}},g_k$, and 
$f_p\not=\overline{g_{q+1}},g_q$.  Then
$$\beta\alpha^{-1}=(e_1\cdots e_j)(g_1\cdots g_k)^{-1} (g_1\cdots 
g_q) (f_1\cdots f_p)^{-1}.$$
Suppose first that $q>k$. Then we have
$$\eqalign{\beta\alpha^{-1}&=e_1\cdots e_j g_{k+1}\cdots g_q 
\overline{f_p}\cdots\overline{f_1}\cr
&=(e_1\cdots e_j g_{k+1}\cdots g_q)\cdot(f_1\cdots f_p)^{-1}\cr
\alpha y&=f_1\cdots f_p g_{q+1}g_{q+2}\cdots\cr
&=(f_1\cdots f_p)\cdot(g_{q+1}g_{q+2}\cdots),}$$
and hence $(\beta\alpha^{-1},\alpha y)$ has as standard form
$$(e_1\cdots e_j g_{k+1}\cdots g_q,\ f_1\cdots f_p,\ 
g_{q+1}g_{q+2}\cdots).$$
It follows that
$$\eqalign{c(\beta\alpha^{-1},\alpha y)&=j+(q-k)-p\cr
&=(j-k)-(p-q)\cr&=c(\beta,y)-c(\alpha,y).\cr}$$
Next, if $q<k$, we may apply the above argument to 
$c\bigl((\alpha,y) (\beta,y)^{-1}\bigr)$. Finally, if $q=k$, we have
$$\eqalign{\beta\alpha^{-1}&=e_1\cdots e_j \overline{f_p}\cdots 
\overline{f_1}\cr
\alpha y&=f_1\cdots f_p g_{k+1}g_{k+2}\cdots.\cr}$$
Let $m$ be the number of cancellations in the product $(e_1\cdots 
e_j)\cdot(\overline{f_p}\cdots\overline{f_1})$.  Then
$$\eqalign{c(\beta\alpha^{-1},\alpha y)&=(j-m)-(p-m)\cr
&=j-p\cr &=(j-k)-(p-q)\cr &=c(\beta,y)-c(\alpha,y).\qed\cr}$$

\smallskip\noindent
{\bf Remark 2.26.}\enspace  Note that $c^{-1}(0)$ is an equivalence relation.  
For 
if $(\beta_1,\beta_2,x)$ is the standard form of the element 
$(\beta_1\beta_2^{-1},\beta_2 x)\in c^{-1}(0)$, then
$$\eqalign{r(\beta_1\beta_2^{-1},\beta_2 x) 
=s(\beta_1\beta_2^{-1},\beta_2 x)&\Longleftrightarrow 
\beta_1x=\beta_2x\cr
&\Longleftrightarrow \beta_1=\beta_2,
\hbox{\ (since\ }\ell(\beta_1)=\ell(\beta_2)),\cr
&\Longleftrightarrow (\beta_1\beta_2^{-1},\beta_2 x) = 
\bigl(\lambda(x),x\bigr)\in c^{-1}(0)^0.\qed\cr}$$

\theorem Proposition 2.27.  Let $K\subseteq 
A(\tildee)\hat$ be a closed $G$-invariant subset, and let
$c:G*K\to\z$ denote the restriction of the
canonical cocycle to $G*K$.  Then
$c^{-1}(0)$ is an AF equivalence 
relation.\endtheorem

\pf For $F\subseteq E$ a finite subgraph, and $n\in\n$, let
$$\displaylines{\quad
H(F,n)=\bigl\{(\alpha,y)\in c^{-1}(0)\bigm| \beta_1,\beta_2\in 
F^*,\ \ell(\beta_1)=\ell(\beta_2)\le n,\ \underline x\in 
E^{**}_+,\hfill\cr
\hskip2true in\hbox{
(where\ $(\beta_1,\beta_2,x)$\ is\ the\ 
standard\ form\ of\ $(\alpha,y)$)}
\bigr\}.\hfill\cr}$$  
Then $H(F,n)$ is a compact-open subequivalence relation of 
$c^{-1}(0)$.  (Compactness follows from the requirement that 
$\underline x$ consist of {\it directed} edges.)
If $F\subseteq F^\prime$ and $n\le n^\prime$, then 
$H(F,n)\subseteq H(F^\prime,n^\prime)$, and
$$\bigcup_{F,n}H(F,n)=c^{-1}(0).$$

We will show that $H(F,n)$ is an elementary groupoid (in the sense of 
\bib\renault\ III.1.1).  Since $c^{-1}(0)^0$ is totally 
disconnected, it will 
follow that $c^{-1}(0)$ is an AF groupoid.  (In \bib\renault\ an AF 
groupoid is defined to be the inductive limit of a {\it sequence} of 
elementary groupoids.  However, since each $H(F,n)$ has totally 
disconnected unit space, $C^*\bigl(H(F,n)\bigr)$ is an AF-algebra.  
Thus $C^*\bigl(c^{-1}(0)\bigr)$ is AF.)

\smallskip
For $y\in H(F,n)^0$, let $M(F,n,y)$ denote the set of ordered pairs 
$(\beta_1,\beta_2)\in F^*\times F^*$ such that
\smallskip
\item{$\bullet$} $\exists x\in K$ with 
$(\beta_1,\beta_2,x)$ equal to the standard form of an element of 
$H(F,n)$,
\item{$\bullet$} $y=\beta_2 x$.
\smallskip\noindent
It is clear that $M(F,n,y)$ is finite.  We claim that there is a 
(compact-open)
neighborhood $U$ of $y$ in $H(F,n)^0$ such that for $y^\prime\in U$, 
$M(F,n,y^\prime)=M(F,n,y)$.  To see this, let $\underline 
y=f_1f_2\cdots$.  If $\ell(\underline y)>n$, let $U=[V(f_1\cdots 
f_{n+1})]\cap H(F,n)^0$.  Then for any $y^\prime\in U$, 
$\underline {y^\prime}=f_1\cdots f_{n+1}\cdots$, and the conclusion 
follows.  If $\ell(y)\le n$, then $\underline y\in \tildee^0$.
Let $R=\Delta_{1,\tildef}(y)$.  This is a finite subset 
of $\Delta_{1,\tildef}(y)$ since $F$ is a finite graph.  We may now take 
$U=[V(y;R)]\cap H(F,n)$.  Then for any $y^\prime\in U$, we must have 
$y^\prime=f_1^\prime\cdots$, where $f_i^\prime=f_i$ for $i\le\ell(y)$, 
and $f_{\ell(y)+1}^\prime\not\in F^1$.  The conclusion now follows 
because $M(F,n,y)\subseteq F^*\times F^*$.

\def\betas{(\beta_1,\beta_2)}
\def\betasprime{(\beta_1^\prime,\beta_2^\prime)}

\smallskip
Now let $y\in H(F,n)^0$ and $U$ be as above.  For $\betas\in 
M(F,n,y)$, let $U\betas=\{\beta_1\beta_2^{-1}\}\times U$.  Note 
that for any $(\alpha,z)\in H(F,n)$ with $z\in U$, there exists 
$\betas\in M(F,n,y)$ with $\alpha=\beta_1\beta_2^{-1}$, and thus such 
that $(\alpha,z)\in U\betas$.  Further, we note that for any 
$\betas\in M(F,n,y)$,
\smallskip
\item{$\bullet$} $U\betas$ is a compact-open $G$-set in $H(F,n)$,
\item{$\bullet$} $s\bigl(U\betas\bigr)=U$,
\item{$\bullet$} For $\betas,\betasprime\in M(F,n,y)$ distinct, 
$r\bigl(U\betas\bigr) 
\cap r\bigl(U\betasprime\bigr)=\emptyset$.
\smallskip\noindent
It follows that
$$\displaylines{\quad
H(F,n)\restrictedto{\bigcup\bigl\{r\bigl(U\betas\bigr) \bigm|
\betas\in M(F,n,y)\bigr\}}
\cong M(F,n,y)^2\times\,U,\hfill\cr}$$
and that $H(F,n)$ is the disjoint union of finitely many such 
restrictions.
Thus $H(F,n)$ is an elementary groupoid.\qed

\smallskip\noindent
\theorem Corollary 2.28.  $C^*(G*K)$ 
is 
nuclear, and coincides with 
$C^*_r(G*K)$.\endtheorem

\pf Let $D=G*K$.  We obtain a 
circle action on $C^*(D)$ in the usual way:
$$\alpha_z(f)(\zeta)=z^{c(\zeta)}f(\zeta),\ f\in C_c(D).$$
(This extends to an automorphism of $C^*(D)$ by Corollary 2.15.)  
Letting 
$H=c^{-1}(0)$, the inclusion $C_c(H)\subseteq C_c(D)\subseteq C^*(D)$ 
extends to an injective $*$-homomorphism $C^*(H)\subseteq C^*(D)$ 
(injectivity follows easily since $C^*(H)$ is an AF algebra).  By 
approximating from within $C_c(H)$, it is easy to see that 
$C^*(H)=C^*(D)^\alpha$, the fixed-point algebra.  Since $C^*(H)$ is 
nuclear, it follows that $C^*(D)$ is nuclear (see \bib\quigg\ for a 
more general result in the context of coactions).

\smallskip
It is easy to check that the above formula for $\alpha_z$ defines a
circle action on $C^*_r(D)$.  Since $C^*(H)=C^*_r(H)$, it follows
that $C^*(H)\subseteq C^*_r(D)$.
Hence if $E$ is the conditional expectation of $C^*(D)$
onto $C^*(H)$, and $\lambda$ is the canonical map of $C^*(D)$
onto $C^*_r(D)$, then $E\circ\lambda=E$.
It now follows that the canonical 
map  is injective.  (Alternatively, one may 
appeal to \bib\anantharamanrenault\ 6.2.14.ii and 6.1.7.) \qed

\smallskip
We now require two simple facts about inclusions of $r$-discrete 
groupoids.

\theorem Lemma 2.29. (i) Let $K$ be an $r$-discrete groupoid and let 
$H$ be an open subgroupoid.  Then $C_c(H)\subseteq C_c(K)$ extends to 
an injective $*$-homomorphism $C^*_r(H)\to C^*_r(K)$.
\smallskip\noindent
(ii) Let $\pi:X\to Y$ be a continuous proper surjection of locally 
compact Hausdorff spaces.  Let $G$ be a discrete groupoid acting 
equivariantly on $X$ and $Y$.  Then $\pi^*:C_c(Y)\to C_c(X)$ extends 
to an injective $*$-homomorphism $C^*_r(G*Y)\to 
C^*_r(G*X)$.\endtheorem

\pf (i) Let $u\in G^0$, and let $\sigma_u$ denote the regular 
representation of $C_c(G*X)$ induced from the point mass at $u$.  
Then $\sigma_u$ acts on $\ell^2(G_u)$ by the formula
$$\sigma_u(f)\xi(\alpha)=\sum_{\beta\in G^{r(\alpha)}}f(\beta) 
\xi(\beta^{-1}\alpha),\ f\in C_c(G),\ \xi \in\ell^2(G_u).$$
We have a partition of $G_u$ into sets invariant for left 
multiplication by $H$:
$$\bigl\{G_u^{G^0\setminus H^0}\bigr\}\cup H\backslash G_u^{H^0}.$$
Choosing a cross-section $F\subseteq G_u^{H^0}$ for
$H\backslash G_u^{H^0}$, we have
$$\sigma_u\restrictedto{C_c(H)}\cong0\oplus\bigoplus_{x\in F} 
\tau_{r(x)},$$
where $\tau_v$ denotes the regular representation of $C_c(H)$ induced 
from the point mass at $v\in H^0$.  Thus for $f\in C_c(H)$,
$$\eqalign{\Vert f\Vert_{C_r^*(G)}&=
\sup_{u\in G^0}\Vert\sigma_u(f)\Vert\cr
&\le\sup_{v\in H^0}\Vert\tau_v(f)\Vert\cr
&=\Vert f\Vert_{C_r^*(H)}.\cr}$$
But if $v\in H^0$, then $\tau_v$ is a subrepresentation of 
$\sigma_v\restrictedto{C_c(H)}$.   Hence $\Vert f\Vert_{C_r^*(G)}
=\Vert f\Vert_{C_r^*(H)}$.
\smallskip\noindent
(ii) Let $\rho:X\to G^0$ and $\lambda:Y\to G^0$ be the fiber maps.  
Let $y\in Y$, let $x\in\pi^{-1}(y)$, and let $\sigma$ (respectively 
$\widetilde\sigma$) denote the regular representation of $C_c(G*Y)$ 
(respectively $C_c(G*X)$) induced from the point mass at $y$ 
(respectively $x$).  Then $\sigma$ (respectively $\widetilde\sigma$) 
acts on $\ell^2(G_{\lambda(y)})$ (respectively 
$\ell^2(G_{\rho(x)})$).  However $\rho(x)=\lambda\pi(x)=\lambda(y)$, 
and it is easy to see that $\widetilde\sigma\circ\pi^*=\sigma$.  Thus 
for $f\in C_c(G*Y)$,
$$\eqalign{\Vert f\Vert_{C_r^*(G*Y)}&=
\sup_{y\in Y}\Vert\sigma_y(f)\Vert\cr
&=\sup_{y\in Y}\Vert\widetilde\sigma_x\circ\pi^*(f)\Vert,
\hbox{\ any\ }x\in\pi^{-1}(y),\cr
&=\sup_{x\in X}\Vert\widetilde\sigma_x\circ\pi^*(f)\Vert\cr
&=\Vert\pi^*(f)\Vert_{C_r^*(G*X)}.\qed\cr}$$

\smallskip\noindent
\theorem Corollary 2.30.  Let $E$ be a directed graph, and let 
$S\subseteq\Sigma(E)$.  Then 
$C^*\bigl(G*\partial(\tildee,\tildes)\bigr)$ and $\t\oh(E,S)$ are 
strongly Morita equivalent.\endtheorem

\pf $\t\g(E,S)$ is an open subgroupoid of 
$G*\partial(\tildee,\tildes)$.  By Lemma 2.29 (i), 
$$C^*_r\bigl(\t\g(E,S)\bigr)\subseteq 
C^*_r\bigl(G*\partial(\tildee,\tildes)\bigr).$$
  But by Corollary 2.28, 
  the full and reduced \cstar-algebras coincide for these 
groupoids.  Therefore
$$\t\oh(E,S)\subseteq C^*\bigl(G*\partial(\tildee,\tildes)\bigr).$$
$$\displaylines{\hbox{Let\ } p=\Chi_{X(E,S)}\in 
M\Bigl(C^*\bigl(G*\partial(\tildee,\tildes)\bigr)\Bigr).
\quad\hbox{Since}\hfill\cr}$$
$$p\; C_c\bigl(G*\partial(\tildee,\tildes)\bigr)p=
C_c\bigl(\t\g(E,S)\bigr),$$
we have by continuity
$$p\; C^*\bigl(G*\partial(\tildee,\tildes)\bigr)p=
\t\oh(E,S).$$
Moreover, $p$ is a full projection since $X(E,S)$ is a transversal in 
$G*\partial(\tildee,\tildes)$.\qed

\smallskip\noindent
{\bf Remark 2.31.}\enspace   If the graph $E$ is countable, then all groupoids 
under discussion are second countable.  In this case the previous 
result follows immediately from \bib\mrw.  Our aim, in proving Corollary 
2.30 as we have done, was not only to establish the result for arbitrary 
graphs, but also to show that its proof is more elementary than would 
be implied by the use of \bib\mrw.

\theorem Lemma 2.32. Let $E$ be a directed graph and let 
$S\subseteq\Sigma$.  Then $C^*\bigl(G*(\Sigma\setminus S)\bigr)$
is isomorphic to 
$\bigoplus\bigl\{\k(H_u)\bigm|u\in\Sigma\setminus S\bigr\}$,
where $H_u=\ell^2\bigl(E^*\cap t\inv(u)\bigr)$.\endtheorem

\pf Since $\Sigma\setminus S$ is a discrete set, we have that 
$G*(\Sigma\setminus S)$ is the disjoint union of the (discrete)
transitive 
groupoids $G*\bigl(E^*\cap t\inv(u)\bigr)$, $u\in\Sigma\setminus S$.  The lemma 
now follows from the fact that the \cstar-algebra of a discrete 
transitive groupoid is elementary. \qed

\theorem Corollary 2.33.  Let $E$ be a directed graph and let 
$S\subseteq \Sigma$.  There is an exact sequence
$$0\longrightarrow I_S\longrightarrow \t\oh(E,S) \longrightarrow 
\oh(E)\longrightarrow0,$$
where $I_S=\bigoplus\bigl\{\k_u\bigm|u\in\Sigma\setminus S\bigr\}$, 
and $\k_u$ is isomorphic to the algebra of compact operators on a 
Hilbert space of dimension $\#\bigl(E^*_+\cap t\inv(u)\bigr)$. 
\endtheorem

\pf This follows from Corollary 2.30, Lemma 2.32 and \bib\renault, Proposition 
II.4.5.(i). \qed

\smallskip
Now let $E_1\subseteq E_2$ be an inclusion 
of directed graphs.  We have 
an inclusion of bundles of directed trees $\tildee_1\subseteq 
\tildee_2$, hence an injective $*$-homomorphism $A(\tildee_1) \to 
A(\tildee_2)$ (Remark 1.19).  Let $S_2\subseteq\Sigma(E_2)$.
We define $S_1\subseteq\Sigma(E_1)$ by
$$S_1=\bigl\{v\in E_1^0\cap S_2\bigm| \Delta_{1,E_2}(v)\subseteq 
E^1_+
\bigr\}.\leqno(*)$$
Then $\widetilde{S_1}$ is the set  defined in 
Lemma 1.20 by the inclusion $\tildee_1\subseteq \tildee_2$ and the set 
$\tildes_2\subseteq\Sigma(\tildee_2)$.  Thus 
the following diagram commutes, has surjective columns 
and  injective rows:
$$\matrix{A(\tildee_1)&\longrightarrow&A(\tildee_2)\cr\cr
\Bigm\downarrow&&\Bigm\downarrow\cr\cr
A(\tildee_1,\tildes_1)&\longrightarrow&A(\tildee_2,\tildes_2).\cr}$$
From Lemma 2.29 we obtain a composition of injective $*$-homomorphisms:
$$C^*\bigl(G(E_1)*\partial(\tildee_1,\tildes_1)\bigr) \to
C^*\bigl(G(E_1)*\partial(\tildee_2,\tildes_2)\bigr) \to
C^*\bigl(G(E_2)*\partial(\tildee_2,\tildes_2)\bigr).$$
We obtain a corresponding injective $*$-homomorphism of Toeplitz graph 
algebras, which we state as a theorem.

\theorem Theorem 2.34.  Let $E_1\subseteq E_2$ be an inclusion of 
directed graphs, and let $S_2\subseteq\Sigma(E_2)$.  Let 
$S_1\subseteq\Sigma(E_1)$ be defined as in $(*)$.  Then there is an 
injective $*$-homorphism
$$\t\oh(E_1,S_1)\to\t\oh(E_2,S_2). \qed$$  \endtheorem

\theorem Theorem 2.35.  Let $E$ be a directed graph.  
Let $\w$ be any collection 
of subgraphs which is directed by inclusion and for which
$$\eqalign{\bigcup_{F\in\w}F^0&=E^0\cr
\bigcup_{F\in\w}F^1&=E^1.\cr}$$
For $S\subseteq\Sigma(E)$, and $F\in\w$, let 
$S_F\subseteq\Sigma(F)$ be as above.  Then
$$\t\oh(E,S)=\lim_{\displaystyle\longrightarrow\atop
\displaystyle\w}\t\oh(F,S_F).$$
\endtheorem

\pf The coherence of the system of $*$-homomorphisms follows from 
Remark 1.21.  Note that \lower5true pt\hbox{$\displaystyle\lim \atop 
{\scriptstyle\longrightarrow\atop 
\scriptstyle\w}$}$A(\tildef,\tildes_F)$ 
is dense in $A(\tildee,\tildes)$ by the Stone-Weierstrass 
theorem, since each $B\in\a(\tildee)$ is the image of a set in 
$\a(\tildef)$ for all large enough $F$.  Since $G(E)=\bigcup_\w 
G(F)$, the equality stated in the theorem holds.\qed

\bigskip\noindent
\centerline{\bf 3.  Structure of Graph $C^*$-algebras.}
\medskip
The results in this section are parallel with similar results in 
other work (e.g. 
{\bf[\exellaca,}\discretionary{}{{\bf\kprr,}}{{\bf\kprr,}}{\bf\kpr]}).
We show here that they follow 
easily in our general framework.
\smallskip\noindent
{\bf Definition 3.1.}\enspace Let $E$ be a directed graph.  A {\it 
cycle} in $E$ is a path $\alpha=e_1\cdots e_n\in E^*_+$ such that 
$o(\alpha)=t(\alpha)$ and $o(e_i)\ne o(e_j)$ for $1\le i<j\le n$.  
An {\it exit} of $\alpha$ is an edge in $\Delta_1\bigl(o(e_i)\bigr) 
\setminus\{e_i\}$ 
for some $i$.  The
cycle $\alpha$ is {\it terminal} if  $\alpha$ 
has no exit. The cycle $\alpha$ is {\it transitory} if it is not 
terminal, and if for each exit $f$ of $\alpha$
 we have $o(e_j)\notin V\oftf$ for all $j$; (in other 
words, $\alpha$ has an exit, but no exit from $\alpha$ leads 
back to $\alpha$).

\theorem Lemma 3.2. Let $E$ be a directed graph.  Let $u\in E^0$, and
  suppose that $E$
has a non-terminal cycle (whose vertices are contained) in $V(u)$.  Then 
the restriction of $G(E)*\partial\tildee$ to
$\bigl[V_\tildee(u)\bigr]\subdel$
contains an open $G$-set $Z$ with $s(\overline Z)\subseteq 
\bigl[V_\tildee(u)\bigr]\subdel$ and $r(\overline 
Z)\subsetnoteq s(Z)$.\endtheorem

\pf Let $\alpha=e_1\cdots e_n$ be a non-terminal cycle with vertices 
in $V(u)$.  
We may assume that $o(\alpha)$ 
emits an exit from $\alpha$.  Let 
$U=\bigl[V_\tildee\bigl(o(\alpha)\bigr)\bigr]\subdel$ and let 
$Z=\{\alpha\}*U$.  Then $Z$ is the required a $G$-set.  (Any directed path 
that begins with an exit from $\alpha$ represents a point of $U$ that 
is not in the range of $Z$.) \qed

\theorem Theorem 3.3.  Let $E$ be a directed graph.  Then $\oh(E)$ is 
an AF algebra if and only if $E$ has no cycles.\endtheorem

\pf Suppose first that $E$ has a cycle $\alpha$.  If $\alpha$ is 
terminal, we may let $N$ be the set of vertices in $\alpha$ and take 
$F$ to be trivial.  Then $(N,F)\in\i(E)$, and the corresponding ideal 
of $\oh(E)$ is isomorphic to a matrix algebra over the continuous 
functions on the circle.  Since $\oh(E)$ has a non-AF ideal, $\oh(E)$ 
is not AF.  If $\alpha$ is not terminal, then let $Z$ be as in
  Lemma 3.2.  Then $\Chi_Z$ is a partial isometry in $\oh(E)$ whose 
  final projection is a proper subprojection of its initial 
  projection.  It follows
that $\oh(E)$ is not AF.

Now suppose that $E$ has no cycles.  Then if $F$ is any finite 
subgraph of $E$, $F^*$ is a finite set, and hence $\t\oh(F,S_F)$ is 
finite dimensional.  Then Theorem 2.35 implies that $\oh(E)$ is AF. \qed

\medskip
We recall from \bib\anantharaman\ that an $r$-discrete groupoid is called {\it locally 
contractive} if for every nonempty open subset $U$ of the unit space there is 
an open $G$-set $Z$ with $s(\overline Z)\subseteq U$ and $r(\overline 
Z)\subsetnoteq s(Z)$ (see also \bib\lacaspielberg).

\theorem Theorem 3.4.  Let $E$ be a directed graph.  Then 
$G(E)*\partial\tildee$ is 
locally contractive if and only if there are no terminal cycles, and 
$V(u)$ contains a cycle for every $u\in E^0$.\endtheorem

\pf We first suppose that $E$ has no terminal cycles and that $V(u)$ 
contains a cycle for every $u\in E^0$.  Then Theorem 3.2 applies to 
every $u\in E^0$.  Since every
nonempty open subset of $\partial\tildee$ contains $V(u)$ for some 
$u\in E^0$, local contractivity follows from Theorem 3.2. 

We now prove the converse.  If $\alpha$ is a terminal cycle, then 
$[\alpha^\infty]$ is an isolated point of $\partial\tildee$, and 
hence is an open set with no contracting subset.  If $u\in E^0$ is 
such that $V(u)$ does not contain a cycle, consider the subgraph $D$ 
of $E$ with $D^0=V(u)$ and $D^1=o\inv(D^0)\cap t\inv(D^0)$.  By 
Theorem 3.3, $\oh(D)$ is AF.  Therefore the restriction of 
$G(E)*\partial\tildee$ to (the open set)
 $\bigl[V_\tildee(u)\bigr]\subdel$ does not 
contain a $G$-set as in the definition of local contractivity.  Hence 
$G(E)*\partial\tildee$ is not locally contractive.
\qed

\smallskip\noindent
{\bf Definition 3.5.}\enspace Let $E$ be a directed graph.  We call 
$E$ {\it cofinal} if the following two properties hold.
\item{(i)} For every $u\in E^0$ and for every $p=e_1e_2\cdots\in 
E^\infty_+$, there is $\alpha\in E^*_+$ with $o(\alpha)=u$ and 
$t(\alpha)=o(e_j)$ for some $j$.
\item{(ii)} $E^0\setminus\Sigma(E)\subseteq V(u)$ for all $u\in E^0$.
\smallskip\noindent
{\bf Remark 3.6.}\enspace The first condition is the notion of cofinality 
(see e.g. \bib\kprr) used in the case of row-finite graphs.  The second 
condition is necessary when there are vertices with infinite exit 
valence.

\theorem Theorem 3.7.  Let $E$ be a directed graph.  Then 
$G(E)*\partial\tildee$ is minimal if and only if $E$ is cofinal. 
\endtheorem

\pf Suppose first that $E$ is cofinal.
 Let $U\subseteq\partial\tildee$ be a nonempty open invariant set. 
We will show that $U=\partial\tildee$. Let $x\in\partial\tildee$. By 
invariance of $U$ we may assume that $x$ is represented by a directed 
path $p\in E^{**}_+$, where $t(p)\in E^0\setminus\Sigma(E)$ if 
$\ell(p)<\infty$. Again by invariance, we may find $u\in E^0$ with 
$\bigl[V_\tildee(u)\bigr]\subdel\subseteq U$.  By cofinality there is 
$\alpha\in E^*_+$ such that $o(\alpha)=u$,
 and $t(\alpha)$ is one of the 
vertices of $p$.  But then $\alpha\inv x\in U$, so that $x\in U$ by 
invariance.

Now suppose that $G(E)*\partial\tildee$ is minimal.  Let $u\in E^0$.  
Put $N_0=V(u)$.  Having defined $N_i$ for $i<i_0$, let 
$$N_{i_0}=\bigl\{v\in \Sigma(E)\bigm| 
t\bigl(\Delta_1(v)\bigr)\subseteq \bigcup_{i<i_0}N_i\bigr\}.$$
Let $N=\cup_i N_i$ and let $F_v=\emptyset$ for $v\in N$.  Refering to 
Definition 1.23, all conditions are evidently true, so that $(N,F)\in 
\i(E)$.  Since $N\ne\emptyset$, minimality and Theorem 2.8 imply that 
$N=E^0$.  Since $N\setminus V(u)\subseteq\Sigma(E)$, condition (ii) of 
cofinality holds.

It is easily shown by induction that if $v\in N$ and if 
$e_1e_2\cdots\in E^\infty_+$ with $o(e_1)=v$, then $t(e_n)\in V(u)$ 
for all $n$ large enough.  Now let $p\in E^\infty_+$.  By minimality 
and Theorem 2.8 there exist $v\in N$, and $\alpha\in E^*$ with 
$t(\alpha)=v$, such that 
$[p]\in\alpha\bigl[V_\tildee(v)\bigr]\subdel$.  Hence there is 
$q=e_1e_2\cdots\in E^\infty_+$ with $o(q)=v$ such that $p=\alpha q$.  
Since $t(e_n)\in V(u)$ for large $n$, condition (i) of cofinality 
holds. \qed

\medskip
We recall that an $r$-discrete groupoid is called {\it 
essentially free} if the set of units having 
trivial isotropy is dense in the space of units.

\theorem Theorem 3.8.  Let $E$ be a directed graph.  Then 
$G(E)*\partial\tildee$ is essentially free if and only if $E$ has no 
terminal cycles.\endtheorem

Before proving Theorem 3.8 we need some lemmas.  We will call a {\it 
circuit} any directed path of positive length whose origin and terminus 
coincide.

\theorem Lemma 3.9.  Let $E$ be a directed graph and let 
$x\in\partial\tildee$.  Then $x$ has nontrivial isotropy if and only 
if $x=[\beta\gamma^\infty]$ for some circuit $\gamma$.\endtheorem

\pf If $x$ has the indicated form, then $\beta\gamma\beta\inv$ 
fixes $x$.  Conversely, suppose $\alpha x=x$ with $\ell(\alpha)>0$.  
Then necessarily $\underline x$ has infinite length.  Let 
$(\alpha,\underline x)\in G(E)*\partial\tildee$ have standard form 
$(\beta_1,\beta_2,y)$.  If $\ell(\beta_1)=\ell(\beta_2)$, then since 
$\beta_1y=x=\beta_2y$ we must have $\beta_1=\beta_2$, contradicting 
the assumption on $\alpha=\beta_1\beta_2\inv$. So we may assume 
$\ell(\beta_1)<\ell(\beta_2)$.  But then from $\beta_1y=\beta_2y$ we 
find that $\beta_2=\beta_1\gamma$.  Then $y=\gamma y$, so 
$y=\gamma^\infty$.  It follows that $x=\beta_1\gamma^\infty$. \qed

\theorem Lemma 3.10.  Let $E$ be a directed graph and let $\alpha\in 
E^*_+$ be a circuit.  Then $\alpha=\beta\gamma\delta$ where $\beta$, 
$\gamma$, $\delta\in E^*_+$ and $\gamma$ is a cycle.\endtheorem

\pf  Let $\gamma$ be the portion of $\alpha$ between two consecutive 
occurrences of a repeated vertex in $\alpha$.  \qed

\theorem Lemma 3.11.  Let $E$ be a directed graph having no terminal 
cycles.  Then for every $u\in E^0$ there is $p\in E^{**}_+$ such 
that $o(p)=u$ and $[p]\in\partial\tildee$ has trivial isotropy. 
\endtheorem

\pf Suppose first that $V(u)$ contains a non-transitory cycle.  Then 
this cycle contains a vertex which emits an edge that is not in the 
cycle, but leads back to a vertex in the cycle.  Thus there are 
$\alpha$, $\beta$, $\gamma$, $\delta\in E^*_+$ such that

\item{(i)} $\beta\gamma$ is the chosen cycle,
\item{(ii)} $o(\delta)=o(\beta$, $t(\delta)=t(\beta)$, and the first 
edge of $\delta$ is not the first edge of $\beta$,
\item{(iii)} $o(\alpha)=u$ and $t(\alpha)=o(\beta)$.

\noindent
Then the path $p=\alpha\beta\gamma\delta\gamma(\beta\gamma)^2 
\delta\gamma(\beta\gamma)^3\cdots$ has the required properties.

Now suppose that $u$ leads only to transitory cycles.  In this case 
there exists a path $p\in E^{**}_+$ with $o(p)=u$, and $t(p)\in 
E^0\setminus\Sigma(E)$ if $\ell(p)<\infty$, such that $p$ contains no 
cycles.  Such a path can be constructed inductively by exiting any 
cycle encountered.  By Lemma 3.10, such a path will contain no 
circuits.  By Lemma 3.9, it will have the required properties. \qed

\medskip\noindent
{\it Proof of Theorem 3.8\/}.\enspace Suppose first that $E$ has no 
terminal cycles.  If $x\in\partial\tildee$ has nontrivial isotropy, 
Lemma 3.9 implies that $x=[\beta\gamma^\infty]$ for some circuit 
$\gamma$.  Let $u=o(\gamma)$.  Let $p$ be as in Lemma 3.11.  Then 
$x_n=[\beta\gamma^np]$ are points in $\partial\tildee$ with trivial 
isotropy that converge to $x$.

Conversely, suppose that $E$ has a terminal cycle, $\alpha$.  Then 
$[\alpha^\infty]$ is an isolated point in $\partial\tildee$, and has 
non-trivial isotropy. \qed

\theorem Corollary 3.12.  Let $E$ be a directed graph, and let 
$S\in\Sigma(E)$.  Then $G(E)*\partial(\tildee,\widetilde S)$ 
is essentially free if and only if 
$E$ has no terminal cycles.\endtheorem

\pf Recall $\Sigma(\tildee)$ is a discrete   clopen 
subset of $A(\tildee)\hat$.  Moreover, every point of 
$\Sigma(\tildee)$ has trivial isotropy.  Therefore 
$\partial(\tildee, \widetilde S) = \partial\tildee\cup 
\bigl(\Sigma(\tildee)\setminus\widetilde S\bigr)$ 
has a dense set of points with 
trivial isotropy if and only if the same is true of 
$\partial\tildee$. \qed

We next recall from \bib\renault, Definition II.4.3, that an $r$-discrete 
groupoid is called {\it essentially principal} if its restriction to 
each closed invariant subset of the unit space is essentially free.
The characterization in graph terms rests on the following lemma.

\theorem Lemma 3.13.  Let $E$ be a directed graph.  Then there exists 
$(N,F)\in\i(E)$ such that $E(N,F)$ has a terminal cycle if and only 
if $E$ has a terminal or transitory cycle.\endtheorem

\pf Suppose that $E$ has a transitory cycle, $\alpha$, and let $C$ 
denote the vertices of the edges in $\alpha$.  Let 
$D=\bigl\{e\in E^1_+\bigm|o(e)\in C,\ t(e)\notin C\bigr\}$.  Then for 
every $u\in t(D)$ we have $V(u)\cap C=\emptyset$.  So we may let 
$N=\cup_{e\in D}V\ofte$, and $F_u=\emptyset$ for $u\in N$.  Then it 
is easy to see that $(N,F)\in\i(E)$, and $\alpha$ becomes a terminal 
cycle in $E(N,F)$.

Conversely, suppose $(N,F)\in\i(E)$ is such that $E(N,F)$ has a 
terminal cycle.  Note that by Definition 1.23 (iia),
$$\bigcup_{u\in N\setminus R(N,F)}V(u)\subseteq N\setminus R(N,F),$$
and that $N\setminus R(N,F)=E^0\setminus E(N,F)^0$.  Now let $\alpha$ 
be a terminal cycle in $E(N,F)$.  With $C$ and $D$ as before, we must 
have
$$\bigcup_{e\in D}V(u)\subseteq E^0\setminus E(N,F)^0.$$
Therefore $\alpha$ is either terminal or transitory in $E$. \qed 

\theorem Theorem 3.14.  Let $E$ be a directed graph.  Then $\g(E)$ is 
essentially principal if and only if $E$ has no terminal or transitory 
cycles.\endtheorem

\pf By Theorem 2.8 and Proposition 1.30, the closed invariant subsets 
of $\partial\tildee$ are of the form 
$$\partial\bigl(E(N,F)\tild,S(N,F)\tild\bigr),$$ for $(N,F)\in\i(E)$.  
By Corollary 3.12, the restriction to such a set is essentially free if 
and only if $E(N,F)$ has no terminal cycles.  By Lemma 3.13, no 
$E(N,F)$ has a terminal cycle if and only if $E$ has no terminal or 
transitory cycles.  \qed

We may assemble some of the above results as follows.

\theorem Theorem 3.15.  Let $E$ be a directed graph.
\item{(i)} $\oh(E)$ is a nuclear $C^*$-algebra.
\item{(ii)} If $E$ has no terminal or transitory cycles, then the 
lattice of ideals in $\oh(E)$ is isomorphic to the lattice of 
invariants $\i(E)$.
\item{(iii)} $\oh(E)$ is simple if and only if $E$ is cofinal and has 
no terminal cycles.
\item{(iv)} $\oh(E)$ is simple and purely infinite if and only if $E$ 
is cofinal, has no terminal cycles, and if $V(u)$ contains a cycle 
for every $u\in E^0$.
\endtheorem

\pf  Statement (i)  follows from Corollaries 2.28 and 2.30.
Statement (ii) follows from Theorem 3.14 and \bib\renault, Theorem II.4.5(iii).
Statement   (iii) follows from (ii) and Theorems 3.7 and 3.8.
Statement (iv) follows from (ii),  Theorem 3.4, and 
\bib\anantharaman\ (or \bib\lacaspielberg, 
whose argument adapts immediately to $r$-discrete groupoids).
\qed
\bigskip
\noindent
{\bf References}
\medskip
\item{\bib\anantharaman} 
C. Anantharaman-Delaroche, Purely infinite \cstar-algebras 
arising from dynamical systems, (preprint).
\item{\bib\anantharamanrenault} 
C. AnantharamanDelaroche and J. Renault, Amenable 
groupoids, (preprint).
\item{\bib\archboldspielberg} 
R. Archbold and J. Spielberg, Topologically free actions 
and ideals in discrete dynamical systems, {\it Proc. Edinburgh Math. 
Soc.\/} {\bf 37} (1993), 119-124.
\item{\bib\bprs} T. Bates, D. Pask, I. Raeburn and W. Szymanski, 
The \cstar-algebras of row-finite graphs, (preprint).
\item{\bib\brenken} B. Brenken, \cstar-algebras of infinite graphs and 
Cuntz-Krieger algebras, (preprint).
\item{\bib\cuntz} J. Cuntz, A class of \cstar-algebras and topological Markov 
chains II:  Reducible chains and the Ext-functor for \cstar-algebras, 
{\it Invent. Math.\/} {\bf 63} (1981), 25-40.
\item{\bib\cuntzkrieger} 
J. Cuntz and W. Krieger, A class of \cstar-algebras and 
topological Markov chains, {\it Invent. Math.\/} {\bf 56} (1980), 
251-268.
\item{\bib\drinentomforde} D. Drinen and M. Tomforde, The \cstar-algebras of 
arbitrary graphs, preprint.
\item{\bib\exellaca} R. Exel and M. Laca, Cuntz-Krieger algebras for infinite 
matrices, {\it J. reine angew. Math.\/} {\bf 512} (1999), 119-172.
\item{\bib\flr} N. Fowler, M. Laca and I. Raeburn,The \cstar-algebras of 
infinite graphs, {\it Proc. Amer. Math. Soc.\/} {\bf 8} (2000), 
2319-2327.
\item{\bib\higgins} P.J. Higgins, The fundamental groupoid of a graph of groups,
{\it J. London Math. Soc. (2)\/} {\bf 13} (1976), 145-149.
\item{\bib\kelley} J.L. Kelley,  {\it General Topology\/} van Nostrand, New 
York, 1955.
\item{\bib\kprr} A. Kumjian, D. Pask, I. Raeburn and J. Renault, 
Graphs, groupoids and Cuntz-Krieger algebras, {\it J. Funct. Anal.\/} {\bf } (), .
\item{\bib\kpr} A. Kumjian, D. Pask and I. Raeburn, 
Cuntz-Krieger algebras of directed graphs, {\it Pacific J. Math.\/} {\bf 
184} (1998), 161-174.
\item{\bib\kumjianpask}  A. Kumjian and D. Pask, \cstar-algebras of directed 
graphs and group actions, {\it Ergodic Theory \& Dyn. Sys.\/} {\bf 19} 
(1999), 1503-1519.
\item{\bib\lacaspielberg} M. Laca and J. Spielberg, Purely infinite \cstar-algebras 
from boundary actions of discrete groups, {\it J. reine angew. Math.\/} 
{\bf 480} (1996), 125-139.
\item{\bib\mrw} P. Muhly, J. Renault and D. Williams, Equivalence and 
isomorphism for groupoid \cstar-algebras, {\it J. Operator Theory\/}
 {\bf 17} (1987), 3-22.
\item{\bib\quigg} J. Quigg, Discrete \cstar-coactions and \cstar-algebraic 
bundles, {\it J. Austral. Math. Soc. (Series A)\/} {\bf 60} (1996), 
204-221.
\item{\bib\quiggsieben} J. Quigg and N. Sieben, \cstar-actions of $r$-discrete 
groupoids and inverse semigroups, (preprint).
\item{\bib\raeburnszymanski} I. Raeburn and W. Szymanski, 
Cuntz-Krieger algebras of infinite graphs and matrices, preprint.
\item{\bib\renault} J. Renault,  {\it A Groupoid Approach to 
\cstar-algebras\/},
Lecture Notes in Mathematics {\bf 793}, Springer-Verlag, Berlin, 1980.
\item{\bib\serre} J. Serre,  {\it Trees\/}, Springer-Verlag, Berlin, 1980.
\item{\bib\spielberg} J. Spielberg, Semiprojectivity for certain purely 
infinite \cstar-algebras, (preprint, available at the Front for the 
Mathematics ArXiv,
{\it http://front.math.ucdavis.edu}
\discretionary
{{\it /math.OA}}
{\it /}
{{\it /math.OA/}}
{\it 0102229\/}).
\item{\bib\szymanskizhang} W. Szymanski and S. Zhang, Infinite simple \cstar-algebras 
and reduced crossed products of abelian \cstar-algebras and free groups, 
{\it Manuscripta Math.\/} {\bf 92} (1997), 487-514.
\bye